\documentclass[a4paper,leqno,11pt]{article}
\usepackage{amsmath,amsthm,amsfonts,amssymb}
\ifx\pdftexversion\undefined
  \usepackage[dvips]{graphicx,color}
\else
  \usepackage[pdftex]{graphicx,color}
\fi
\usepackage[francais]{babel}
\usepackage[latin1]{inputenc}%

\catcode`\×=\active
\def×{\csname times\endcsname} 

\newcounter{annexe} 

\makeatletter 
\newcommand\annexe{%
    \@startsection{annexe}{1}{\z@}%
    {10mm plus 6mm minus 3mm}{\baselineskip}%
    {\normalfont\normalsize\scshape\centering}}
\@addtoreset{subsection}{annexe}
\newcommand{\annexemark}[1]{}
\renewcommand{\@seccntformat}[1]{\if e\if x\if e\if\if a#1\fi\fi\fi\fi Annexe \csname the#1\endcsname\quad%
\else\csname the#1\endcsname\quad\fi}
\renewcommand\theenumi{\@roman\c@enumi}
\renewcommand\theenumii{\@alph\c@enumii}

 \makeatother

\newenvironment{enumebis}{\begin{itemize}\item[]\begin{enumerate}}{\end{enumerate}\end{itemize}}

\usepackage{pifont} \renewcommand{\qed}{\flushright\ding{111}} 

\let\ifanglais\iffalse

\def\and{ et }
\def\andname{et}
\def\EMdash{\leavevmode\hbox to 7.5mm{\vrule height .63ex depth -.59ex
    width 5.4mm\hfill}}
\def\email#1{{\it e-mail :\/} #1}
\newcommand{\norm}[1]{\| #1 \|}
\long\def\address#1{\removelastskip\nobreak\vskip\baselineskip\noindent
\vbox\bgroup\def\\{\vskip 1.0mm\parindent=.75cm\relax}\small
\parindent=0pt\obeylines {#1}}
\def\endaddress{\egroup}
\edef\deuxpoints{:}
\makeatletter
\renewenvironment{proof}[1][\!]{\par
  \pushQED{\qed}%
  \normalfont \topsep6\p@\@plus6\p@\relax
  \trivlist
  \item[\hskip\labelsep
        \itshape
    Démonstration #1\@addpunct{.}]
\item[]\ignorespaces
}{%
  \begin{flushright}\popQED\end{flushright}\endtrivlist\@endpefalse
}   
\makeatother
\renewcommand{\:}{\colon}

\def\R{{\mathbb R}}  
 \def\H{{\mathbb H}}

\newcommand{\e}{\varepsilon}
\newcommand{\ed}[1]{\textrm{d} #1}

\def \ignorepar{\afterassignment\ignoreparaux \let\next=} \def
\ignoreparaux{\ifx\next\par \let\next\ignorepar \fi \next}

 \def\vol{{\rm vol}}
\def\argth{{\rm Argth}}

\newcommand{\volumename}{volume} 
\newcommand{\ofname}{de}

\newcommand{\Inname}{Dans} 
\newcommand{\pagesname}{pages}
\renewcommand{\leq}{\leqslant} \renewcommand{\geq}{\geqslant} 

\newtheoremstyle{mesthm}
  {10pt plus 1pt minus 1pt}
  {9pt minus 6pt}
  {\slshape}
  {0.5cm}
  {\bfseries}
  {.}
  {1ex}
  {}
\newtheoremstyle{mesdefi}
  {6pt plus 1pt minus 1pt}
  {6pt plus 1pt minus 1pt}
  {}
  {0.5cm}
  {\bfseries}
  {.}
  {1ex}
  {}%
  
\theoremstyle{mesthm}
\newtheorem{lema}{\ifanglais{\large L}emma\else{\large L}emme\fi}
\newtheorem{theo}[lema]{\ifanglais{\large T}heorem\else {\large
    T}héorème\fi}

\newtheorem{prop}[lema]{{\large P}roposition}  \newtheorem{cor}{{\large C}orollaire}[lema]

\newtheorem{lemme}{\ifanglais{\large L}emma\else{\large L}emme\fi}

\theoremstyle{mesdefi}
 \newtheorem{ex}[lema]{\ifanglais{\large
    E}xample\else{\large E}xemple\fi}

\makeatletter
\newenvironment{rem}{\vskip 10pt plus 1pt minus 1pt \noindent\begin{minipage}{\textwidth}%
  \hskip 0.5cm{\bf\ifanglais{\large R}emark\else{\large R}emarque\fi.}\hskip 1ex
\ignorepar}%
{\end{minipage}\vskip 9pt minus 6pt}
\newenvironment{thsn}[1]{\protected@edef\@currentlabel%
  {\csname p@cor\endcsname#1}
\vskip 10pt plus 1pt minus 1pt  \noindent\begin{minipage}{\textwidth}%
  \hskip 0.5cm{\bf\ifanglais{\large T}heorem \else{\large T}héorème
    \fi \textbf{#1}.}\hskip 1ex
\itshape \ignorepar}%
{\end{minipage}\vskip 9pt minus 6pt}
  
  \title{L'AIRE DES TRIANGLES IDÉAUX \\ EN GÉOMÉTRIE DE HILBERT}
  
  \author
{B.~\sc{Colbois}, C.~\sc{Vernicos}\footnote{Partiellement financé par le projet
européen ACR OFES numéro 00.0349 et la bourse FNRS 20-65060.01}\ \and  P.~\sc{Verovic} }

\begin{document} 

\maketitle

\begin{abstract}
  L'objet de cet article est l'étude de l'aire des triangles idéaux pour la
  géométrie de Hilbert d'un domaine convexe de $\R^n$.
  Les résultats que nous obtenons donnent d'une part une caractérisation de la géométrie
  hyperbolique dans l'ensemble des géométries de Hilbert,
  et d'autre part une minoration optimale, indépendante du convexe,
  de l'aire de Hilbert des triangles idéaux qui caractérise les domaines triangulaires du plan.
  En outre, sous certaines conditions géométriques, nous établissons une majoration de cette aire
  dont nous montrons qu'elle doit dépendre du convexe.
\end{abstract}

\section*{Introduction}

  Le concept de simplexe idéal joue un rôle important dans l'étude des
var­iétés riemanniennes à courbure négative.  Par exemple, J.~Barge et
É.~Ghys obtiennent la caractérisation suivante de la géométrie
hyperbolique plane comme conséquence de leur résultat sur la
cohomologie bornée (voir \cite{bg}, p.~511)~:

\begin{theo} \label{theobg}
  Soit $g$ une métrique riemannienne de courbure négative ou nulle sur
  une surface $S$ compacte,
  connexe et orientable. \\
  Si les triangles idéaux du revêtement universel de $S$ ont tous la
  même aire, alors $(S,g)$ est de courbure constante.
\end{theo}

\noindent
Signalons que pour une surface riemannienne complète et simplement
connexe à courbure négative ou nulle dont tous les triangles idéaux ont une aire finie,
on ne sait toujours pas s'il existe un analogue de ce résultat.

Dans la première partie du présent travail, nous obtenons une
caractérisation de la géométrie hyperbolique parmi les géométries de
Hilbert en terme d'aire des triangles idéaux (voir le théorème
\ref{theo1} ci-dessous).  Cette caractérisation peut être considérée
comme une généralisation du théorème précédent dans un cadre quelque
peu différent. Puis nous étudions les problèmes de minoration et
majoration de l'aire des triangles idéaux.

\begin{figure}[htpb] 
  \centering \input{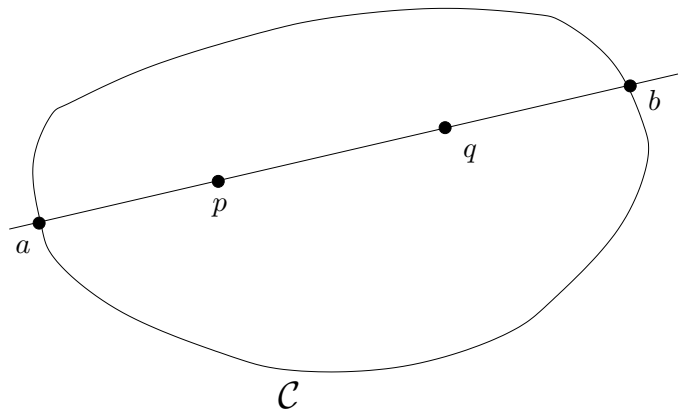}
  \caption{Distance de Hilbert \label{dintro}}
\end{figure}

Avant d'énoncer précisément nos résultats, rappelons qu'une géométrie
de Hilbert $(\mathcal C, d_{\mathcal C})$ est la donnée d'un ouvert
non vide, convexe et borné $\mathcal C$ de $\mathbb R^n$ --- que nous
appelerons \textit{domaine convexe} --- muni de la distance de Hilbert
$d_{\mathcal C}$ définie de la manière suivante~: pour tous points
distincts $p$ et $q$ dans $\mathcal{C}$, la droite passant par $p$ et
$q$ rencontre le bord $\partial \mathcal C$ de $\mathcal C$ en deux points
$a$ et $b$ tels que $p$ soit entre $a$ et $q$ et $q$ soit entre $p$ et
$b$ (figure~\ref{dintro}).  On définit alors
$$
d_{\mathcal C}(p,q) = \frac{1}{2} \ln [a,p,q,b],
$$
où $[a,p,q,b]$ est le birapport de $(a,p,q,b)$, c'est-à-dire
$$
[a,p,q,b] = \frac{\| q-a \|}{\| p-a \|} × \frac{\| p-b \|}{\| q-b
  \|} > 1,
$$
en désignant par $\| \cdot \|$ la norme euclidienne canonique sur
$\mathbb R^n$.  On pose également $d_{\mathcal C}(p,p) = 0$ (voir \cite{dhilbert}, appendice I).  $\\$

Remarquons tout de suite que si $\mathcal{C}$ et $\mathcal{C}'$ sont
deux domaines convexes de $\R^{n}$ tels que leurs images
respectives $\widehat{\mathcal{C}}$ et $\widehat{\mathcal{C}'}$ dans
l'espace projectif ${\mathbb P}^{n}(\R)$ vérifient
$\widehat{\mathcal{C}'} = A(\widehat{\mathcal{C}})$, où $A$ est une
homographie de ${\mathbb P}^{n}(\R)$ --- donc conserve le
birapport de quatre points de ${\mathbb P}^{n}(\R)$---, alors les
géométries de Hilbert $(\mathcal{C} , d_{\mathcal{C}})$ et
$(\mathcal{C}' , d_{\mathcal{C}'})$ sont isométriques.

Dans toute géométrie de Hilbert $(\mathcal C, d_{\mathcal C})$, le
segment de droite reliant deux points quelconques du convexe $\mathcal
C$ est un segment géodésique pour $d_{\mathcal C}$ (au sens de
\cite{bridsonh},~p.~4) et $(\mathcal C, d_{\mathcal C})$ est un espace
métrique géodésique dont la topologie est celle induite par la
topologie canonique de $\R^{n}$.  Ceci dit, en général, le
segment reliant deux points n'est pas l'unique géodésique entre
ceux-ci, cette unicité étant néanmoins satisfaite lorsque le bord $\partial
\mathcal C$ de $\mathcal C$ est une hypersurface de classe $C^2$ dans
$\R^{n}$ dont la courbure de Gauss est partout non nulle --- on
dira alors que $\mathcal C$ est un \textit{convexe strict}.  Notons
enfin que cette condition d'être un convexe strict n'est pas
nécessaire pour avoir unicité du segment géodésique --- voir une
discussion détaillée de ce point dans \cite{so},~§~1.2.2.

Par ailleurs, on peut mettre sur tout domaine convexe $\mathcal{C} \subset\R^n$
une métrique de Finsler $C^0$, notée $F_{\mathcal C}$, en
procédant comme suit~: si $p \in \mathcal C$ et $v \in T_{p}\mathcal C =\R^n$
avec $v \neq 0$, la droite passant par $p$ et dirigée par
$v$ coupe $\partial \mathcal C$ en deux points $p_{\mathcal C}^{+}$ et
$p_{\mathcal C}^{-}$~; on pose alors
$$
F_{\mathcal C}(p,v) = \frac{1}{2} \| v \| \biggl(\frac{1}{\| p -
  p_{\mathcal C}^{-} \|} + \frac{1}{\| p - p_{\mathcal C}^{+}
  \|}\biggr) \quad \textrm{et} \quad F_{\mathcal C}(p , 0) = 0.
$$
Cette métrique de Finsler est liée à la distance de Hilbert
$d_{\mathcal C}$ par le fait que
$$
F_{\mathcal C}(p,v) = \frac{\ed{}}{\ed{t}}\Big{|}_{ t = 0} d_{\mathcal
  C}(p , p + tv)
$$
et que
$$
d_{\mathcal C}(p,q) = \inf \biggl\{ \int_{0}^{1} \! F_{\mathcal C}(\sigma(t)
, \sigma'(t)) \, \ed{t} ~\Big{|}~ \sigma \in \Omega^{1}({\mathcal C} , p , q) \biggr\},
$$
où
$$
\Omega^{1}({\mathcal C} , p , q) = \bigl\{ \sigma \: [0,1] \longrightarrow {\mathcal C}
~\big{|}~ \sigma \text{ de classe } C^{1} \text{ avec } \sigma(0) = p \text{
  et } \sigma(1) = q \bigr\}\text{.}
$$
Grâce à cette métrique de Finsler, on construit une mesure
borélienne $\mu_{\mathcal C}$ sur $\mathcal C$ (qui correspond en fait
à la mesure de Hausdorff de l'espace métrique $(\mathcal C,
d_{\mathcal C})$ --- voir \cite{bbi}, exemple~5.5.13 ) que nous allons
expliciter.

Pour chaque $p \in \mathcal C$, soient $B_{\mathcal C}(p) = \{v \in \R^n ~|~ F_{\mathcal C}(p,v) < 1 \}$
la boule unité ouverte
de $T_{p}\mathcal C = \R^n$ pour la norme $F_{\mathcal
  C}(p,\cdot)$ et $\omega_{n}$ le volume euclidien de la boule unité
ouverte de l'espace euclidien canonique $\R^n$.  En considérant
la fonction (densité) $h : \mathcal C \longrightarrow \R$ donnée par $h(p)
= \omega_{n}/\vol\bigl(B_{\mathcal C}(p)\bigr),$ où $\vol$ est la mesure de
Lebesgue canonique sur $\R^n$, on définit $\mu_{\mathcal C}$ ---
que nous appelerons \textit{mesure de Hilbert} sur $\mathcal{C}$ ---
par
$$
\mu_{\mathcal C}(A) = \int_{A} h(p) \ed{\vol(p)}
$$
pour tout borélien $A$ de $\mathcal C$.

Lorsque $\mathcal{C}$ est un ellipsoïde, $(\mathcal C,
d_{\mathcal C})$ correspond au modèle projectif (ou modèle de Klein)
de la géométrie hyperbolique, et on peut penser aux géométries de
Hilbert $(\mathcal C, d_{\mathcal C})$ comme à une généralisation
naturelle de l'espace hyperbolique.  Une question commune à de
nombreux travaux récents (voir \cite{so}, \cite{so2}, \cite{be},
\cite{cv}, \cite{kn} et leurs références) est de déterminer les
propriétés de l'espace hyperbolique dont héritent les géométries de
Hilbert et de trouver des caractérisations de l'espace hyperbolique
parmi celles-ci.

Le premier résultat de cet article est l'obtention d'une telle
caractérisation grâce à l'aire de Hilbert des triangles idéaux. À
cause de la non unicité des géodésiques pour $d_{\mathcal C}$ entre
deux points d'un domaine convexe $\mathcal C \subset \R^{n}$, un
triangle de $(\mathcal C , d_{\mathcal C})$ ne peut être défini à
l'aide des segments géodésiques de $d_{\mathcal C}$ qui joignent ses
sommets.  C'est pourquoi nous convenons de définir tout d'abord un
triangle $T = a b c$ de $\R^n$ comme l'intérieur de l'enveloppe convexe
affine ouverte de trois points non alignés $a, b, c \in \R^{n}$.  Un
tel triangle sera alors un triangle de $(\mathcal C , d_{\mathcal C})$
si ses sommets sont dans $\mathcal C$ et un triangle idéal de
$(\mathcal C , d_{\mathcal C})$ si ses sommets sont dans $\partial \mathcal
C$ et s'il est inclus dans~$\mathcal C$.

Dans le cas d'un convexe strict, cela équivaut à la définition usuelle
d'un triangle idéal d'un espace métrique uniquement géodésique, en
particulier de l'espace hyperbolique $\H^n$ dans lequel il est connu que tous les
triangles idéaux sont isométriques avec une aire (hyperbolique)
commune égale à $\pi$.  En fait, nous allons montrer que cette
propriété de l'aire caractérise $\H^n$ parmi les
géométries de Hilbert de $\R^n$ :

\begin{theo} \label{theo1}
  Étant donné une géométrie de Hilbert $(\mathcal{C} , d_{\mathcal
    C})$ avec $\mathcal C\in \R^n$, on a~:
  \begin{enumebis}
  \item Tous les triangles idéaux de $(\mathcal C,d_{\mathcal C})$
    sont d'aire constante si, et seulement si, $\mathcal{C}$ est un
    ellipsoïde --- auquel cas cette aire constante vaut~$\pi$.

  \item Si $\mathcal C$ n'est pas un ellipsoïde, il existe des triangles
    idéaux de $(\mathcal C,d_{\mathcal C})$ d'aire strictement plus
    grande que $\pi$ et d'autres d'aire strictement plus petite que
    $\pi$.
    \end{enumebis}
\end{theo}

\begin{rem}
Ici, et dans toute la suite de ce travail,
l'aire d'un triangle (idéal ou pas) de $(\mathcal C,d_{\mathcal C})$
est son aire pour la mesure de Hilbert de $(\mathcal C \cap P , d_{\mathcal C \cap P})$,
où $P$ est l'unique plan vectoriel de $\R^n$ contenant le triangle.  
\end{rem}

La démonstration du théorème \ref{theo1}, donnée dans la première
partie de cet article, est simple et purement géométrique.

Dans la seconde partie, nous obtenons une minoration uniforme de
l'aire des triangles idéaux avec caractérisation du cas d'égalité :

\begin{theo} \label{theolb}
  Étant donné une géométrie de Hilbert $(\mathcal{C} ,
  d_{\mathcal{C}})$ avec $\mathcal C \subset \R^{n}$, on a~:

  \begin{enumebis}
  \item L'aire de tout triangle idéal de $(\mathcal{C} ,
    d_{\mathcal{C}})$ est au moins égale à $\pi^3\!/24$.

  \item Si $n = 2$ et s'il existe un triangle idéal de $(\mathcal{C} ,
    d_{\mathcal{C}})$ d'aire égal à $\pi^3\!/24$, alors $\mathcal{C}$ est
    un domaine triangulaire.
  \end{enumebis}
\end{theo}

Remarquons que le cas d'égalité caractérise bien la géométrie
de $\mathcal{C}$, puisque tous les domaines
triangulaires du plan munis de leurs géométries de Hilbert sont isométriques.

Enfin, dans la troisième partie, nous montrons que la recherche d'une majoration de l'aire des
triangles idéaux donne lieu à une situation différente et
plus contrastée.  En effet, le corollaire~\ref{coin} ci-dessous
fournit des géométries de Hilbert qui possèdent des triangles idéaux
d'aire infinie, de sorte qu'il est illusoire de chercher un majorant
de l'aire des triangles idéaux commun à \textit{toutes} les géométries
de Hilbert à l'instar du théorème~\ref{theolb}.  L'exemple~\ref{exhuit}
montre également que cette impossibilité persiste même en
se restreignant à l'ensemble des convexes stricts de $\R^n$.

Cependant, lorsqu'on considère un convexe strict \textit{fixé}
$\mathcal{C}$ de $\R^n$, nous prouvons qu'il existe néanmoins un
majorant (dépendant de $\mathcal{C}$) de l'aire de \textit{tous} les
triangles idéaux de $(\mathcal C , d_{\mathcal C})$~:

\begin{theo} \label{upbound}
  Soit $\mathcal{C}$ un convexe strict de $\R^n$.  Alors il existe une
  constante $\alpha = \alpha(\mathcal{C}) > 0$ telle que tout triangle idéal
  de $(\mathcal{C} , d_{\mathcal{C}})$ a une aire au plus égale à
  $\alpha$.
\end{theo}

\section{Préliminaires}
\subsection{Quelques propriétés élémentaires}

Nous débutons par une liste de faits simples et généraux dont nous
ferons abondamment usage.

\begin{figure}[htpb]
  \centering \input{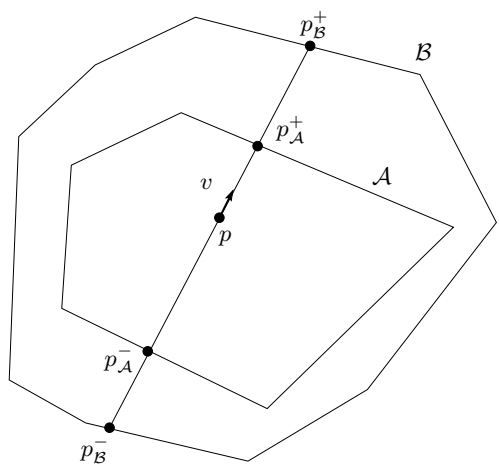}
\caption{Comparaison des distances et mesures de Hilbert
  de deux domaines convexes emboîtés \label{norme}}
\end{figure}

\begin{prop} \label{pptinclusion}
  Soient $({\mathcal A},d_{\mathcal A})$ et $({\mathcal B},
  d_{\mathcal B})$ des géométries de Hilbert telles que ${\cal A} \subset
  {\cal B} \subset \R^{n}$. Alors~:

    \begin{enumebis}
    \item \label{v1} Les métriques de Finsler $F_{\cal A}$ et $F_{\cal
        B}$ vérifient $F_{\mathcal B}(p,v) \leq F_{\mathcal A}(p,v)$
      pour tous $p\in \mathcal A$ et $v \in \mathbb R^{n}$ non nul,
      l'égalité ayant lieu si, et seulement si, $p_\mathcal{A}^{-} =
      p_\mathcal{B}^{-}$ et $p_\mathcal{A}^{+} = p_\mathcal{B}^{+}$
      (figure~\ref{norme}).

    \item \label{v2} Pour tous $p, q \in \mathcal{A}$, on a
      $d_\mathcal{B}(p,q) \leq d_\mathcal{A}(p,q)$.

    \item \label{v3} Pour tout $p \in \mathcal{A}$, on a $\vol \! \left(
        B_\mathcal{A}(p) \right) \leq \vol \! \left( B_\mathcal{B}(p)
      \right)$, avec égalité si, et seulement si, $\mathcal{A} =
      \mathcal{B}$.

    \item \label{v4} Pour tout borélien $A$ de $\mathcal{A}$, on a
      $\mu_{\mathcal{B}}(A) \leq \mu_{\mathcal{A}}(A)$, avec égalité si,
      et seulement si, $\mathcal{A} = \mathcal{B}$.
\end{enumebis}
\end{prop}

\begin{proof}

  Il suffit de prouver l'assertion (\ref{v1}) qui implique toutes les
  autres propriétés. Or, elle découle directement du fait que
  pour tous $p \in \mathcal{A}$ et $v \in \R^n$, $v \neq 0$, on a
  $$
  \| p - p_{\mathcal{A}}^{+} \| \leq \| p - p_{\mathcal{B}}^{+} \|
  \quad \textrm{et} \quad \| p - p_{\mathcal{A}}^{-} \| \leq \| p -
  p_{\mathcal{B}}^{-} \|,
  $$
  l'égalité ayant lieu si, et seulement si, $p_\mathcal{A}^{-} =
  p_\mathcal{B}^{-}$ et $p_\mathcal{A}^{+} = p_\mathcal{B}^{+}$.
  \end{proof}

Grâce à cette proposition, on va pouvoir estimer la mesure de Hilbert
d'un domaine convexe du plan inclus dans un domaine carré, ce dernier
présentant l'avantage d'être suffisamment simple pour permettre
des calculs effectifs.

\subsection{Estimation de l'aire par comparaison avec le domaine carré}

L'estimation de l'aire de Hilbert d'un convexe de $\R^2$ revient à estimer le volume
euclidien de la boule unité ouverte pour la métrique de Finsler en chaque point du convexe.
Lorsque le convexe est un carré, on obtient~:

\begin{prop} \label{squarevol}
  Soit le domaine carré $\mathcal{S} = \{ (x , y) \in \R^{2} ~|~ |x| <
  1 \textrm{ et } |y| < 1 \}$.  Alors pour tout $p = (x , y) \in
  \mathcal{S}$, on a
  $$
  2 (1 - x^{2}) (1 - y^{2}) \leq \vol(B_{\mathcal{S}}(p)) \leq 4 (1 -x^{2}) (1 - y^{2}),
  $$
où $B_{\mathcal{S}}(p)$ est la boule unité ouverte de $T_pS = \R^2$ pour la norme $F_S(p , \cdot)$.
\end{prop}

\begin{proof}

Étant donné $p \in \mathcal{S}$, la preuve consiste à vérifier que
la boule $B_{\mathcal{S}}(p)$ est d'une part incluse dans
un rectangle $\mathcal{R}$ dont les côtés sont parallèlles à ceux du carré $\mathcal{S}$,
et d'autre part contient un losange dont les sommets sont les points de contact entre $\mathcal{R}$
et $B_{\mathcal{S}}(p)$.

  Puisque $\mathcal{S}$ est symétrique par rapport aux axes de
  coordonnées, il suffit de se restreindre à $p \in [0 , 1[ × [0 , 1[$.
  \begin{itemize}
  \item[\textbullet] Soit $v = (a , b) \in \R^{2}$ non nul tel que $|a| \leq\frac{1}{2} (1 - x)$ 
  et $b \geq 0$, de sorte que la demi-droite $p +
  \R_{-} v$ (resp. $p + \R_{+} v$) coupe $\partial
  \mathcal{S}$ sur la droite d'équation $y = -1$ (resp. $y = 1$) en un
  point $p_{\mathcal{S}}^{-}$ (resp. $p_{\mathcal{S}}^{+}$).  Il
  résulte alors du théorème de Thalès que
  $$
  \frac{b}{1 + y} = \frac{\norm{v}}{\norm{p - p_{\mathcal{S}}^{-}}}
  \quad \textrm{ et } \quad \frac{b}{1 - y} = \frac{\norm{v}}{\norm{p
      - p_{\mathcal{S}}^{+}}}~,
  $$
  d'où
  $$
  F_{\mathcal{S}}(p , v) = \frac{1}{2} \! \left( \frac{b}{1 + y} +
    \frac{b}{1 - y} \right) = \frac{b}{1 - y^{2}}~,
  $$
  ce qui donne l'implication $v \in B_{\mathcal{S}}(p) \Longrightarrow b < 1 - y^{2}$.

  Ainsi, $B_{\mathcal{S}}(p)$ étant symétrique par rapport à $0$, on a
  $$
  B_{\mathcal{S}}(p) \cap \left\{  \Bigl[-\frac{1}{2} (1 - x) , \frac{1}{2} (1 -
  x)\Bigr] × \R\right\}  \subset \R × [-(1 - y^{2}) , (1 - y^{2})]
  $$
  et par suite
  $$
  B_{\mathcal{S}}(p) \subset \R × [-(1 - y^{2}) , (1 - y^{2})]
  $$
puisque $B_{\mathcal{S}}(p)$ est convexe.

  De la même façon, on montre que
  $$
  B_{\mathcal{S}}(p) \subset [-(1 - x^{2}) , (1 - x^{2})] × \R.
  $$
  Par conséquent, on obtient
  $$
  B_{\mathcal{S}}(p) \subset [-(1 - x^{2}) , (1 - x^{2})] × [-(1 - y^{2}) , (1 - y^{2})],
  $$
  ce qui entra\^{\i}ne la deuxième inégalité de la
  proposition~\ref{squarevol}.

  \item[\textbullet] On remarque par ailleurs que les points $(1 - x^{2} , 0)$ et
  $(0 , 1 - y^{2})$ sont dans l'adhérence de $B_{\mathcal{S}}(p)$, qui
  est convexe et symétrique par rapport à $0$, d'où il résulte que
  l'enveloppe convexe des points $(1 - x^{2} , 0)$, $(0 , 1 - y^{2})$,
  $-(1 - x^{2} , 0)$ et $-(0 , 1 - y^{2})$ est dans
  $\overline{B_{\mathcal{S}}(p)}$.  Comme le volume euclidien de cette
  enveloppe convexe --- qui est un losange --- est égal à $2 (1 - x^{2}) (1 - y^{2})$, on en
  déduit la première inégalité de la proposition~\ref{squarevol}.

\end{itemize}
\end{proof}

\begin{rem}
À titre indicatif, on peut aisément voir que la boule $B_{\mathcal{S}}(p)$ est un octogone
lorsque $p$ n'est pas sur les diagonales de $\mathcal{S}$,
sinon $B_{\mathcal{S}}(p)$ est un hexagone
si $p \neq 0$ et un carré si $p = 0$.
\end{rem}

De cette estimation, nous pouvons alors tirer deux conséquences
utiles concernant l'aire de Hilbert des triangles idéaux.

\begin{cor} \label{plat}
  Soient $\mathcal{C}$ un domaine convexe du plan tel que $\partial
  \mathcal{C}$ contient un segment ouvert $]a , b[$ et $p \in
  \mathcal{C}$.  Pour chaque $t \in \, ]0 , 1[$, notons $m_{a}(t) = (1
  - t) p + t a$ et $m_{b}(t) = (1 - t) p + t b$.
  Alors, pour $0 < s < t$, si $A(t,s)$
  désigne l'enveloppe convexe des points $m_{a}(t)$, $m_{a}(s)$,
  $m_{b}(t)$ et $m_{b}(s)$, on a $\displaystyle \lim_{t \to 1}
  \mu_{\mathcal{C}}\bigl(A(t,s)\bigr) = + \infty$ lorsque $s$ est fixé.
\end{cor}

\begin{proof}

  Après transformation affine, on se ramène au cas où $p = 0$ et
  $\mathcal{C}$ est inclus dans le carré $\mathcal{S}$ de la
  proposition~\ref{squarevol} avec $a = (-x_{0} , 1)$ et $b = (x_{0} ,
  1)$ pour un certain $x_{0} \in \, ]0 , 1[$ (figure~\ref{square1}).

\begin{figure}[htpb]
  \centering \includegraphics{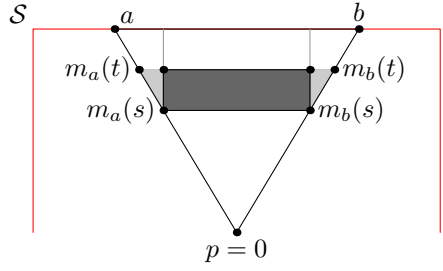}
  \caption{Cas d'un convexe dont le bord contient un segment \label{square1}}
\end{figure}

  Alors, pout tous $s,t \in \, ]0 , 1[$ tels que $s\leq t$, le rectangle de sommets $m_a(s)=(-sx_{0} , s)$,
    $m_b(s)=(s x_{0} , s)$, $m_a(t)=(-sx_{0} , t)$
  et $m_b(t)=(s x_{0} , t)$ est inclus dans $A(s,t)$, d'où il résulte que
  $$
  \mu_{\mathcal{S}}(A(s,t)) \geq 2 \! \int_{0}^{s x_{0}} \! \left(
    \int_{s}^{t} \frac{\pi}{\vol\bigl(B_{\mathcal{S}}(x , y)\bigr)}~\ed{y} \right) \! \ed{x}.
  $$
  Mais, d'après la proposition~\ref{squarevol}, on a
  $\vol(B_{\mathcal{S}}(x , y)) \leq 4 (1 - x^{2}) (1 - y^{2})$ pour
  tout $p = (x , y) \in \mathcal{S}$, ce qui entra\^{\i}ne que
  $$
  \mu_{\mathcal{S}}(A(s,t)) \geq \frac{\pi}{2} \! × \! \left( \int_{0}^{sx_{0}} \!
  \frac{\ed{x}}{1 - x^{2}} \right) \! × \! \left( \int_{s}^{t} \! \frac{\ed{y}}{1 - y^{2}} \right)\!,
  $$
  c'est-à-dire,
\begin{eqnarray*}
\mu_{\mathcal{S}}(A(s,t))
& \geq &
\frac{\pi}{2} \! × \! \textrm{Argth}(s x_{0})
\! × \! \left[ \textrm{Argth}(t) - \textrm{Argth}(s) \right]. \\
\end{eqnarray*}

Par conséquent, en fixant $s$, on obtient $\displaystyle \lim_{t \to 1}
\mu_{\mathcal{S}}\bigl(A(s,t)\bigr) = + \infty$.

Comme $\mathcal{C} \subset \mathcal{S}$, on a finalement $\displaystyle
\lim_{t \to 1} \mu_{\mathcal{C}}\bigl(A(s,t)\bigr) = + \infty$ d'après la
proposition~\ref{pptinclusion}~(\ref{v4}).  \end{proof}

\begin{cor} \label{coin}
  Soient $\mathcal{C}$ un domaine convexe du plan et $\omega \in \partial
  \mathcal{C}$ tels qu'il existe deux droites d'appui distinctes de
  $\mathcal{C}$ en $\omega$.

  Alors, pour tous points distincts $p, q \in \mathcal{C}$, on a
  $\mu_{\mathcal{C}}(p \omega q) = + \infty$, où $p \omega q$ est le triangle de
  sommets $p$, $q$ et $\omega$.
\end{cor}

\begin{proof}

  On va montrer que tout triangle de $\mathcal{C}$ dont un sommet est un \og\,coin\,\fg{}
  de $\mathcal{C}$ peut être pensé comme une demi-bande affine ouverte du plan.

Par transformation affine, on se ramène au cas où $\mathcal{C}$ est
  inclus dans le carré $\mathcal{S}$ de la proposition~\ref{squarevol}
  avec $\omega = (1 , 1)$ et les droites $(\omega p)$ et $(\omega q)$ symétriques
  l'une de l'autre dans la réflexion par rapport à la droite $(0 \omega)$ et tel que
$p_{0} , q_{0} \in \cal{C}$, où $p_{0}$ et $q_{0}$ sont respectivement les points
d'intersection de la droite d'équation $x + y = 1$ avec $(\omega p)$ et $(\omega q)$.

\begin{figure}[htpb]
  \centering \includegraphics{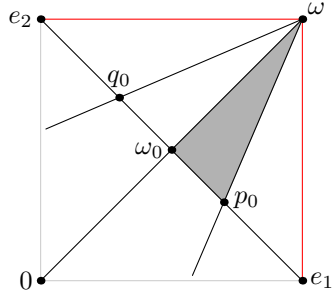}
  \caption{Cas d'un convexe possédant un \og\,coin\,\fg{}}
\end{figure}

En notant $e_{1} = (1 , 0)$, $e_{2} = (0 , 1)$ et $\omega_{0} = (1/2 , 1/2)$,
il existe donc $t_{0} \in \, ]0 , 1[$ tel que
$p_{0} = (1 - t_{0}) \omega_{0} + t_{0} e_{1}$ et
$q_{0} = (1 - t_{0}) \omega_{0} + t_{0} e_{2}$.

Considérons alors
$\Delta = \{ (x , y) \in \R^{2} ~|~ x + y > 1 \textrm{~et~} y < x < 1 \} \subset {\cal S}$
et le $C^{\infty}$-difféomorphisme $f : \Delta \longrightarrow \R^{*}_{+} × \R^{*}_{+}$ défini
par
$$
f(x , y) = (X , Y) = (\argth(t) , \argth(s)) \textrm{,}
$$
où $t , s \in \, ]0 , 1[$ sont tels que
$(x , y) = (1 - s) [(1 - t) \omega_{0} + t e_{1}] + s \omega$.

L'image par $f$ du triangle $\omega_{0} \omega p_{0} \subset \Delta$
est ainsi la bande $]0 , \argth(t_{0})[ × \R^{*}_{+}$
dont on va montrer que l'aire euclidienne usuelle --- qui est infinie ---
est plus petite que l'aire de $\omega_{0} \omega p_{0}$ pour la mesure
de Hilbert $\mu_{\cal S}$.

Un calcul simple donne
$$
t = \frac{y - x}{y + x - 2} \qquad \text{et} \qquad s = x + y - 1,
$$
ce qui entraîne que le jacobien de $f$ en $(x , y) \in \Delta$ vaut
$$
\text{Jac}(f)(x,y)=\frac{1}{2 (x + y) (1 - x) (1 - y)}~\text{.}
$$

En vertu de la deuxième inégalité de la proposition~\ref{squarevol} et du fait que
$$
(1 + x)(1 + y) \leq 3 (x + y) \textrm{~~pour tout~~} (x , y) \in \Delta \textrm{,}
$$
il en résulte que
\begin{multline*}
+\infty
=
\int_{]0 , \argth(t_{0})[ × \R^{*}_{+}} \ed{X} \ed{Y}
=\\
\int_{\omega_{0} \omega p_{0}} \frac{\ed{x} \ed{y}}{2 (x + y) (1 - x) (1 - y)}
\leq
\frac{6}{\pi} \, \mu_{S}({\omega_{0} \omega p_{0}})
\end{multline*}

D'autre part, puisque
$\mu_{S}({\omega_{0} \omega p_{0}}) \leq \mu_{S}({q_{0} \omega p_{0}})$
et que
$\mu_{S}((q \omega p) \backslash (q_{0} \omega p_{0}))$ est finie
--- la partie $(q \omega p) \backslash (q_{0} \omega p_{0})$ étant compacte ---,
on en déduit que $\mu_{S}(q \omega p) = +\infty$.

Enfin, comme $\mathcal{C} \subset \mathcal{S}$, la
proposition~\ref{pptinclusion}~(\ref{v4}) achève la preuve du
corollaire~\ref{coin}.  \end{proof}

\section{Caractérisation de la géométrie hyperbolique par l'aire de Hilbert des triangles idéaux}

En préliminaire à la démonstration du théorème~\ref{theo1},
rappelons le théorème suivant qui est un résultat classique de géométrie convexe
que nous énonçons en dimension deux et dont la preuve se trouve dans
\cite{John} ou \cite{levy}, Lecture~$3$, Theorem~3.1, p.~$13$--$19$.

\begin{theo}[Ellipse de John] \label{johnny}
  Soit $\mathcal{C}$ un domaine convexe du plan.

  Il contient une unique ellipse ouverte d'aire euclidienne maximale,
  l'ellipse de John de $\mathcal{C}$, dont le bord a au moins trois
  points de contact avec $\partial \mathcal{C}$.

  Par dualité, $\mathcal{C}$ est aussi inclus dans une unique ellipse
  ouverte d'aire euclidienne minimale dont le bord a au moins trois
  points de contact avec $\partial \mathcal{C}$.
\end{theo}

Nous allons maintenant donner la preuve du théorème~\ref{theo1} qui
indique comment l'aire de Hilbert des triangles idéaux permet de
caractériser l'espace hyperbolique $\H^n$ parmi
toutes les géométries de Hilbert de $\R^n$.

\begin{proof}[du théorème~\ref{theo1}]

Commençons par faire la preuve lorsque ${\mathcal C} \subset \R^2$.

  Si $\mathcal C$ est une ellipse, l'espace métrique $(\mathcal C ,
  d_{\mathcal C})$ est isométrique au modèle projectif de Klein du
  plan hyperbolique (voir par exemple \cite{bepe},~p.~2) qui a tous
  ses triangles idéaux d'aire égale à~$\pi$.

  Si $\mathcal C$ est n'est pas une ellipse, soit $\mathcal{E}_i$
  l'unique ellipse ouverte d'aire euclidienne maximale incluse dans le
  convexe $\mathcal C$ --- donnée par le théorème \ref{johnny} et
  appelée ellipse de John de $\mathcal C$. L'ellipse $\mathcal{E}_i$ ayant au moins
  trois points de contact avec $\partial \mathcal C$, on peut
  considérer le triangle $T_i$ dont les sommets sont ces trois points
  (figure~\ref{johnsel}).

\begin{figure}[htpb]
  \centering \includegraphics{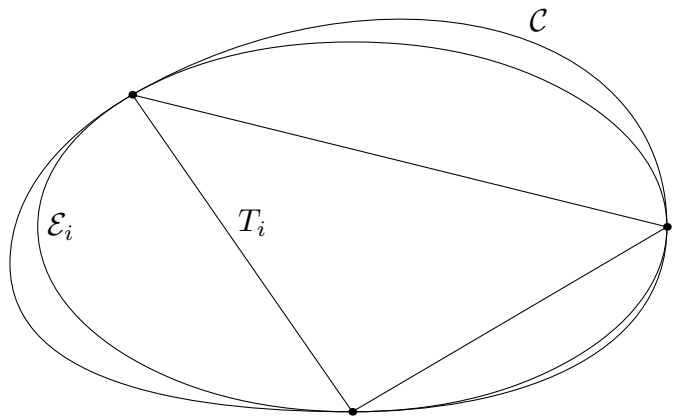}
  \caption{Ellipse de John \label{johnsel}}
\end{figure}

Pour la géométrie de Hilbert associée à l'ellipse de John
$\mathcal{E}_i$, le triangle $T_i$ est idéal et d'aire égale à
\hbox{$\mu_{\mathcal{E}_i}(T_i) = \pi$}.  Par conséquent, comme
$\mathcal{E}_i$ est strictement incluse dans $\mathcal C$, on a
$\mu_{\mathcal{C}}(T_i) < \pi$ en vertu de la
proposition~\ref{pptinclusion}~(\ref{v4}).

D'autre part, considérons l'unique ellipse ouverte $\mathcal{E}_e$
d'aire euclidienne minimale contenant $\mathcal C$ (duale de
$\mathcal{E}_i$). D'après le théorème~\ref{johnny}, son bord possède également
au moins trois points en
commun avec celui de $\mathcal C$, ce qui définit un triangle $T_{e}$.
\begin{enumerate}
\item Si $T_{e}$ est un triangle idéal de $(\mathcal C , d_{\mathcal
  C})$, alors $\mu_{\mathcal{C}}(T_e) > \pi$ puisque $\mathcal{E}_e$
contient strictement $\mathcal C$.

\item Si $T_{e}$ n'est pas un triangle idéal de $(\mathcal{C} ,
d_{\mathcal{C}})$, alors l'un des côtés du triangle $T_{e}$ est inclus
dans $\partial \mathcal{C}$, ce qui implique que l'on peut obtenir un
triangle idéal de $(\mathcal{C} , d_{\mathcal{C}})$ dont l'aire est
arbitrairement grande en vertu du corollaire~\ref{plat}.

Enfin, dans le cas où $\mathcal{C} \subset \R^n$, on fait ce qui précède
dans chaque intersection de $\mathcal{C}$ avec un plan vectoriel de $\R^n$,
sachant que $\mathcal{C}$ est un ellipsoïde si, et seulement si, chacune
de ces intersections est une ellipse.
\end{enumerate}
\end{proof}

\section{Bornes sur l'aire des triangles idéaux en géométrie de Hilbert}

Nous sommes à présent naturellement amenés à nous demander si l'aire
des triangles idéaux d'une géométrie de Hilbert est contrôlée.

\subsection{Du côté de la minoration}

En ce qui concerne la minoration de l'aire des triangles idéaux, nous
avons le résultat global énoncé au théorème~\ref{theolb} qui est
valable pour n'importe quel domaine convexe de $\R^n$.  Pour démontrer
ceci, on va étudier au préalable le cas particulier où le domaine
convexe est un triangle de $\R^2$.

\begin{lema} \label{Letriangle}
  Soit $\Delta$ un domaine convexe triangulaire du plan.

  Alors tous les triangles idéaux de $(\Delta , d_{\Delta})$ ont une aire au
  moins égale à $\pi^3\!/24$ et seul le triangle idéal de sommets les
  milieux des côtés de $\Delta$ a une aire égale à ce minimum.
\end{lema}

Ce lemme étant assez technique, nous ne donnerons que les étapes
de sa preuve, renvoyant à l'annexe~\ref{annexeA1} pour les détails.
Mais auparavant, montrons comment ce lemme implique le théorème~\ref{theolb}.

\begin{proof}[du théorème~\ref{theolb} à l'aide du lemme \ref{Letriangle}]

Considérons d'abord le cas où $\mathcal{C} \subset \R^2$.

  Soient $T = a b c$ un triangle idéal de $(\mathcal{C} ,
  d_{\mathcal{C}})$ et $D_{a}$, $D_{b}$, $D_{c}$ des droites d'appui
  du convexe $\mathcal C$ en $a$, $b$ et $c$ respectivement.
  
   D'un point de vue projectif, il s'agit d'une même et unique situation. Cependant,
d'un point de vue affine --- celui que nous avons suivi jusqu'ici ---,
trois cas se présentent~:

\begin{enumerate}
\item Si $D_{a}$, $D_{b}$ et $D_{c}$ définissent un domaine convexe
  triangulaire $\Delta$ qui contient $\mathcal{C}$, alors on obtient le point (i)
  du théorème~\ref{theolb}
  en appliquant la proposition~\ref{pptinclusion}~(\ref{v4}) (avec
  $\mathcal{A} = \mathcal{C}$ et $\mathcal{B} = \Delta$) et le
  lemme~\ref{Letriangle}.
En outre, si $\mathcal{C}$ est strictement inclus dans $\Delta$,
alors $\mu_{\mathcal{C}}(T) > \mu_{\Delta}(T) \geq \pi^3\!/24$
d'après la proposition~\ref{pptinclusion},
d'où $\mu_{\mathcal{C}}(T) > \pi^3\!/24$, ce qui donne le point (ii)
du théorème~\ref{theolb} en contraposant.

\item Si $D_{a}$ et $D_{b}$ sont parallèles, on plonge le plan affine
contenant le convexe $\mathcal{C}$ dans son complété projectif (voir par exemple \cite{berger}, § 5.1)
identifié naturellement à $\mathbb{P}^{2}(\R)$, dans lequel les deux droites $D_a$ et
$D_b$ se coupent.
En considérant alors une droite $D$
qui est l'image par le plongement d'une parallèle à $D_c$ contenue dans le demi-plan déterminé par
$D_c$ et ne contenant pas $\mathcal{C}$, on est ramené au point (i)
dans le nouveau plan affine $\mathbb{P}^{2}(\R) \!\! \setminus \!\! D$.

\item Si $D_{a}$, $D_{b}$ et $D_{c}$ définissent un domaine convexe
  triangulaire qui ne contient pas $\mathcal{C}$, on peut supposer que le triangle en question
est dans le demi-plan déterminé par $D_a$ et ne contenant pas $\mathcal{C}$.
Dans ce cas, on plonge le plan affine contenant
$\mathcal{C}$ dans $\mathbb{P}^{2}(\R)$.
En considérant alors une droite $D$ image par le plongement d'une droite parallèle à $D_a$
contenue dans le demi-plan déterminé par $D_a$, ne contenant pas $\mathcal{C}$
et qui rencontre le domaine triangulaire, on est ramené au point (i)
dans le nouveau plan affine $\mathbb{P}^{2}(\R) \!\! \setminus \!\! D$.

Enfin, dans le cas où $\mathcal{C} \subset \R^n$, on applique ce qui précède
à chaque intersection de $\mathcal{C}$ avec un plan vectoriel de $\R^n$.
\end{enumerate}

\end{proof}

\begin{proof}[du lemme~\ref{Letriangle}]

  Partant d'un domaine triangulaire $\Delta = m p q \subset \R^{2}$ et
  d'un triangle idéal $T = a b c$ de $(\Delta , d_{\Delta})$, la preuve va se
  faire en trois étapes.

\smallskip
\noindent\textbf{Étape~1~:} 
Elle se résume au lemme suivant, dont la preuve est en annexe~\ref{annexeA1}.
\begin{lema}\label{etape1}
Soit $\Delta_0$ le domaine triangulaire de sommets $0$, $e_1$ et $e_2$,
où $(e_1 , e_2)$ est la base canonique de $\R^2$.
Pour chaque $\alpha \in \, ]0 , 1 / 2]$, notons $T(\alpha)$ le triangle idéal de $(\Delta_0,d_{\Delta_0})$
dont les sommets sont
\begin{align*}
  a(\alpha) &= (\alpha , 1 - \alpha), &b(\alpha)& = (0 , 1 - \alpha) &&\text{et} &c(\alpha) &=
  (\alpha , 0)\text{.}
\end{align*}
Alors il existe $\alpha\in\,]0,1/2]$ et une transformation affine de $\R^2$ qui envoit simultanément
le domaine $\Delta$ sur le domaine $\Delta_0$ et le triangle idéal $T$ de $(\Delta , d_{\Delta})$
sur le triangle idéal $T(\alpha)$ de $(\Delta_0,d_{\Delta_0})$.
\end{lema}

\smallskip
\noindent\textbf{Étape~2~:} Sachant que pour chaque $p = (x , y) \in \Delta_{0}$, la
boule unité ouverte $B_{\Delta_{0}}(p)$ de $T_{p} \Delta_{0} = \mathbb R^2$
pour la norme $F_{\Delta_{0}}(p , \cdot)$ est un hexagone décrit dans
\cite{dlharpe}, p.~106--107, le calcul de
la mesure de Hilbert $\mu_{\Delta_{0}}$ (que nous ne détaillerons pas) nous
donne
$$
\ed{\mu_{\Delta_{0}}(\,p)} =\frac{\pi}{12} × \frac{\ed{x}\,\ed{y}}{xy(1-x-y)}~{\rm.}
$$
L'application $\mathcal{A} : \, ]0 , 1/2] \longrightarrow \R$ définie
par $\mathcal{A}(\alpha) = \mu_{\Delta_{0}}(T(\alpha))$ est strictement
décroissante de sorte que son minimum est atteint en $\alpha = 1/2$
seulement, ce qui correspond au triangle idéal $T(1/2)$ de $(\Delta_{0} ,
d_{\Delta_{0}})$ dont les sommets sont les milieux des côtés de $\Delta_{0}$
--- voir les calculs dans l'annexe \ref{annexeA2}.

\smallskip
\noindent\textbf{Étape~3~:} Montrons que l'aire de Hilbert de $T(1/2)$ est
égale à $\pi^3\!/24$.

D'après l'annexe \ref{annexeA2}, cela revient à calculer
$$
{\mathcal F}(0)= -2\!\int_0^1 \frac{\ln(1-x)}{x}\ed{x} + 2\!\int_0^1
\frac{\ln(1+x)}{x}\ed{x}
$$
puisque $\mathcal{A}(1/2) = \frac{\pi}{12} \mathcal{F}(0)$.

Pour déterminer le premier terme, on développe en série entière par
rapport à $\e \in \, ]0 , 1[$ la quantité
\begin{eqnarray*}
  \label{eq1}
  {\mathcal F}_1(\e) & = & \int_0^{1-\e} \frac{\ed{x}}{x} \ln\Bigl( \frac{1}{1-x} \Bigr) \nonumber \\
     & = & \int_0^{1-\e} \frac{\ed{x}}{x} \sum_{k=1}^{+\infty}\frac{x^k}{k} \nonumber \\
     & = & \sum_{k=1}^{+\infty} \frac{1}{k^2} (1-\e)^k.
\end{eqnarray*}
En faisant $\e \to 0$, le théorème de convergence dominée de Lebesgue
nous permet alors d'obtenir ${\mathcal F}_1(0) = \pi^2 \!/ 6$.

Pour évaluer le second terme de ${\mathcal F}(0)$, on introduit la
fonction
$$
{\mathcal F}_2(\e) = \int_0^{1-\e} \frac{\ed{x}}{x} \ln(1+x)
$$
définie pour $\e \in \, ]0 , 1[$ et qui, également à l'aide d'un
développement en série entière, fournit
$$
{\mathcal F}_2(0) = -\sum_{k=1} ^{+\infty}\frac{(-1)^k}{k^2}~.
$$
En remarquant alors que ${\mathcal F}_1(0) - {\mathcal F}_2(0) = \frac{1}{2} {\mathcal F}_1(0)$, on
obtient ${\mathcal F}_2(0) = \pi^{2} /!/ 12$ et par suite ${\mathcal F}(0) = {\mathcal F}_1(0) +
{\mathcal F}_2(0) = \pi^2 /!/ 2$.  Finalement, l'aire de Hilbert de $T(1/2)$ vaut
$\mathcal{A}(1/2) = \frac{\pi}{12} × \frac{\pi^{2}}{2} =
\frac{\pi^{3}}{24}$~.  \end{proof}

\subsection{Du côté de la majoration}

Pour ce qui est de la majoration de l'aire des triangles idéaux,
remarquons tout d'abord qu'il existe des géométries de Hilbert planes
dans lesquelles on peut trouver des triangles idéaux d'aire aussi
grande que l'on veut, et même d'aire infinie, comme le montrent les
corollaires~\ref{plat} et \ref{coin} de la première partie.

Néanmoins, avec quelques hypothèses de régularité, on évite les
triangles idéaux d'aire infinie~:

\begin{prop} \label{prop:airefinie}
  Soit $\mathcal{C}$ un domaine convexe de $\R^n$ dont le bord est une
  hypersurface de classe $C^2$.
  Alors tout triangle idéal de $(\mathcal{C} , d_{\mathcal{C}})$ a une aire finie.
\end{prop}

\begin{proof}

Considérons d'abord le cas où $\mathcal{C} \subset \R^2$.

  Soit $T = a b c$ un triangle idéal de $(\mathcal{C} ,
  d_{\mathcal{C}})$
  dont on oriente les sommets dans le sens trigonométrique. \\
  Comme $\partial \mathcal C$ est de classe $C^{2}$, il existe $r > 0$ et
  des disques ouverts euclidiens $D(a)$, $D(b)$ et $D(c)$ de rayon $r$
  tangents au bord de $\mathcal C$ en $a$, $b$ et $c$ respectivement
  et inclus dans $\mathcal C$.  En considérant le sommet $a$,
  désignons par $a'$ et $a''$ les points d'intersection du bord de
  $D(a)$ avec les segments $]a , b[$ et $]a , c[$ respectivement.
  Ainsi le triangle $a a' a''$ a ses sommets orientés dans le sens
  trigonométrique et a $\partial D(a)$ pour cercle euclidien circonscrit.
  En procédant de même avec les sommets $b$ et $c$, on obtient les
  triangles $b b' b''$ et $c c' c''$.

\begin{figure}[htbp]
  \centering \includegraphics{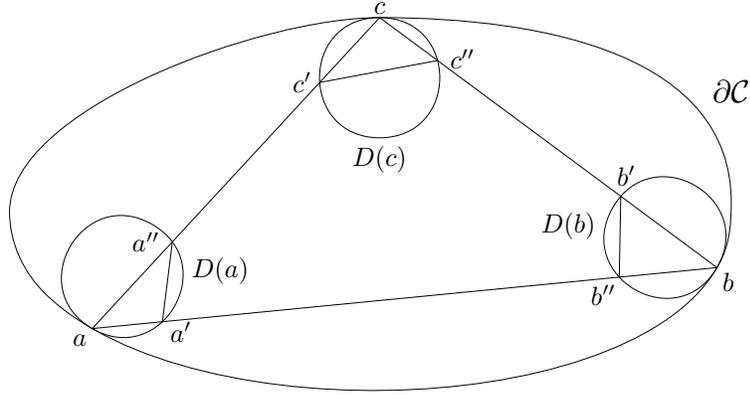}
  \caption{L'aire des triangles idéaux est finie
    dès que $\partial \mathcal{C}$ est $C^{2}$\label{airefinie}}
\end{figure}

L'adhérence du complémentaire de la réunion des triangles $a a' a''$,
$b b' b''$ et $c c' c''$ dans $T = a b c$ est alors un compact inclus
dans $\mathcal{C}$, donc d'aire de Hilbert finie.  En outre, comme $a
a' a''$ est un triangle idéal de $(D(a) , d_{D(a)})$, on a
$\mu_{D(a)}(a a' a'') = \pi$ (géométrie hyperbolique plane) et par suite
$\mu_{\mathcal{C}}(a a' a'') \leq \pi$ en vertu de la
proposition~\ref{pptinclusion}~(\ref{v4}).  Comme il en est de même
avec $b b' b''$ et $c c' c''$, la proposition~\ref{prop:airefinie} en
découle lorsque $\mathcal{C} \subset \R^2$.

Dans le cas où $\mathcal{C} \subset \R^n$, on fait ce qui précède
dans chaque intersection de $\mathcal{C}$ avec un plan vectoriel de $\R^n$.
\end{proof}

Comme on l'a déjà vu au corollaire~\ref{plat}, dès qu'un domaine
convexe $\mathcal{C}$ du plan a un bord --- même de classe $C^2$ --- qui
contient un segment ouvert, alors l'aire des triangles idéaux de
$(\mathcal{C} , d_{\mathcal{C}})$ peut être arbitrairement grande.  En
revanche, lorsque $\mathcal{C}$ est un convexe strict, nous montrons
que ceci ne peut pas se produire~:

\begin{thsn}{\ref{upbound}}
  Soit $\mathcal{C}$ un convexe strict de $\R^n$.  Alors il existe une
  constante $\alpha = \alpha(\mathcal{C}) > 0$ telle que tout triangle idéal
  de $(\mathcal{C} , d_{\mathcal{C}})$ a une aire au plus égale à
  $\alpha$.
\end{thsn}
 
Avant de donner la preuve de ce théorème, notons cependant que la
constante $\alpha(\mathcal{C})$ n'admet pas de majoration uniforme en
$\mathcal{C}$, même dans l'ensemble des convexes stricts du plan,
comme le montre l'exemple suivant.

\begin{ex} \label{exhuit}
  Considérons le carré $\mathcal{S} = \{ (x , y) \in \R^{2} ~|~
  |x| < 1 \textrm{ et } |y| < 1 \}$ ainsi que son homothétique $t
  \mathcal{S}$ avec $t \in \, ]1 / 2 , 1[$ arbitraire.  Si $\mathcal{C}$
  est un convexe du plan tel que $t \mathcal{S} \subset \mathcal{C} \subset
  \mathcal{S}$, désignons par $a_{\mathcal{C}}$, $b_{\mathcal{C}}$ et
  $c_{\mathcal{C}}$ les points d'intersection de $\partial \mathcal{C}$ avec
  les segments fermés joignant $0$ à $a = (-1 , 1)$, $b = (1 , 1)$ et
  $c = (0 , -1)$ respectivement (figure~\ref{entre2}).
  
  Alors $T_{\mathcal{C}} = a_{\mathcal{C}} b_{\mathcal{C}}
  c_{\mathcal{C}}$ est un triangle idéal de $\mathcal{C}$ qui contient
  la partie $A(1 / 2 , t)$ définie au corollaire~\ref{plat} avec $p = 0$ et $s = 1 / 2$.
  Le même corollaire affirmant que $\lim_{t \to 1} \mu_{\mathcal{S}}\bigl(A(1/2,t)\bigr) = +\infty$, il
  en résulte que $\lim_{t \to 1} \mu_{\mathcal{C}}\bigl(A(t)\bigr) = +\infty$ et par
  suite, pour tout entier $n > 0$, il existe $t \in \, ]0 , 1[$ tel que
  tout domaine convexe $\mathcal{C}$ du plan avec $t \mathcal{S} \subset
  \mathcal{C} \subset \mathcal{S}$ (qu'il soit strict ou pas) vérifie
  $\mu_{\mathcal{S}}(T_{\mathcal{C}}) > n$.

\begin{figure}[hbtp]   
  \centering \includegraphics{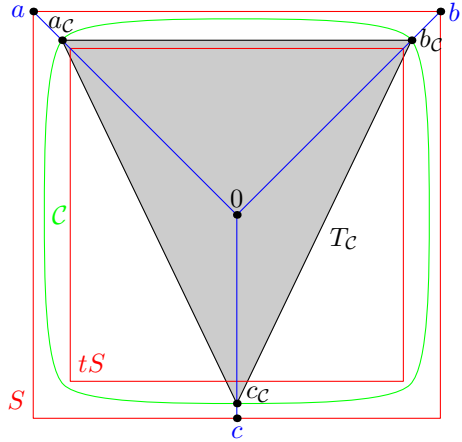}
  \caption{L'aire des triangles idéaux ne peut être majorée
    uniformément en $\mathcal{C}$\label{entre2}}
  \end{figure}
\end{ex}


Afin de démontrer le théorème~\ref{upbound},
pour lequel on se ramène au cas où $\mathcal{C} \subset \R^2$
par intersection avec des plans vectoriels de $\R^n$,
nous allons utiliser la
distance euclidienne canonique $d$ sur $\R^{2}$ et écrire
$\mathcal C$ comme la réunion du compact $K_{\delta} = \{ p \in \mathcal{C}
~|~ d(p , \partial \mathcal{C}) \geq \delta \}$ --- où la constante $\delta = \delta(\mathcal
C) > 0$ sera précisée ultérieurement --- et de son complémentaire $V_{\delta}
= \mathcal C \!\!\setminus\!\! K_{\delta}$.  Pour majorer l'aire d'un triangle idéal
quelconque $T$ de $(\mathcal{C} , d_{\mathcal{C}})$, il suffira alors
de majorer l'aire de la partie du triangle $T$ hors du compact $K_\delta$,
c'est-à-dire $T \cap V_{\delta}$.  Pour cela, nous allons inclure $T \cap
V_{\delta}$ dans la réunion d'un certain nombre $N = N(\mathcal C) > 0$ de
triangles, chacun d'eux étant contenu dans un disque ouvert inclus
dans $\mathcal C$.

La proposition~\ref{pptinclusion}~(\ref{v4}) permet alors de majorer
l'aire de $T \cap V_{\delta}$ par $N\pi$ en comparaison avec la géométrie du
plan hyperbolique (qui est, rappelons-le, la géométrie de Hilbert d'un
disque ouvert).

Finalement, on aura $\mu_{\mathcal C}(T) \leq N\pi + \mu_{\mathcal
  C}(K_{\delta})$.

La preuve de ce résultat reposera sur plusieurs lemmes techniques
démontrés dans l'annexe~\ref{annexeB}.

Tout d'abord, le fait que $\mathcal C$ soit un convexe strict assure
l'existence de deux constantes $r > 0$ et $R > 0$ telles que le cercle
de rayon $2 r$ roule à l'intérieur de $\overline{\mathcal C}$ et que
$\partial \mathcal C$ roule à l'intérieur du disque fermé de rayon $R$ (voir
\cite{blaschke} et \cite{cv}, p.~3).  Cela va nous permettre de ramener une partie de la
preuve du théorème~\ref{upbound} à des considérations sur les cordes
de deux cercles euclidiens embo\^{\i}tés et tangents données au
lemme~\ref{Cordes}.

Puis, à l'aide du lemme~\ref{Distbord} ci-dessous, on étudiera les
cordes de $\partial \mathcal C$ (c'est-à-dire les segments fermés reliants
deux points distincts de $\partial \mathcal C$) en les comparant aux cordes
des cercles euclidiens de rayon $r$ (resp. $R$) tangents
intérieurement (resp. extérieurement) à $\overline{\mathcal C}$.  Lors
de la preuve du théorème~\ref{upbound}, ces cordes seront les côtés
des triangles idéaux de $(\mathcal{C} , d_{\mathcal{C}})$ et la
constante $\delta = \delta(\mathcal C)$ ne dépendra que de $r$ et $R$.

Dans la suite, pour tous $t > 0$ et $\omega \in \partial \mathcal C$, on
désignera par $\Gamma_{t}(\omega)$ le cercle de rayon $t$ tangent à $\partial
\mathcal C$ en $\omega$ et inclus dans le demi-plan fermé de
$\R^{2}$ contenant $\omega$ dans son bord et dans lequel se trouve
$\mathcal C$.  Aussi, le disque ouvert correspondant sera noté
$D_{t}(\omega)$.

\begin{lema} \label{Distbord} 
Soit $\mathcal{C}$ un convexe strict du plan.
  Pour tous points distincts $a$ et $b$ de $\partial \mathcal C$, on a~:

    \begin{enumerate}
    \item Il existe un unique point d'intersection $a'$ entre $]a ,
      b[$ et $\Gamma_{r}(a)$.

    \item La distance euclidienne de $a'$ au bord de $\mathcal C$ est
      minorée en fonction de $d(a , b)$, $r$ et $R$ uniquement~:
      $$
      d(a' , \partial \mathcal C) \geq \frac{r}{4 R^{2}} d(a , b)^{2}.
      $$
    \end{enumerate}
\end{lema}

Les preuves de ce lemme et du lemme~\ref{Cordes} qui l'implique,
seront données dans l'annexe~\ref{annexeB}, tout comme le résultat suivant, qui
fournit la clé du théorème~\ref{upbound}~:

\begin{lema} \label{Recouvrement} 
Soient $\mathcal{C}$ un convexe strict du plan ainsi que $a$ et $b$
deux points distincts de $\partial \mathcal C$ tels que $d(a , b) \leq r$.
Alors, l'unique rectangle ouvert $S(a , b)$ de base le segment $]a, b[$
et de hauteur $r$ inclus dans $\mathcal{C}$ vérifie
$$
\mu_{\mathcal C}(S(a , b)) \leq 2 \pi E\!\!\left( \!\frac{2 R}{r} \!\right)\!\text{.}
$$

        
      
      
      
       
\end{lema}

À partir de maintenant, posons $\delta = \frac{r^3}{4 R^2}$ et
introduisons
$$
V_{\delta} = \{ p \in \mathcal C ~|~ d(p , \partial \mathcal C) < \delta \}
$$
dont le complémentaire $K_{\delta} = {\mathcal C} \backslash V_{\delta}$ est
compact, donc d'aire $\mu_{\mathcal C}(K_{\delta})$ finie.

Quitte à diminuer $r$ et/ou augmenter $R$, on peut supposer que $\delta$
est suffisamment petit pour que $K_{\delta}$ soit convexe en utilisant
l'exponentielle normale de la sous-variété $\partial \mathcal C$ de
$\R^{2}$ muni de sa métrique riemannienne canonique (voir
\cite{cv},~p.~5)

\begin{figure}[hbtp] 
  \centering \includegraphics{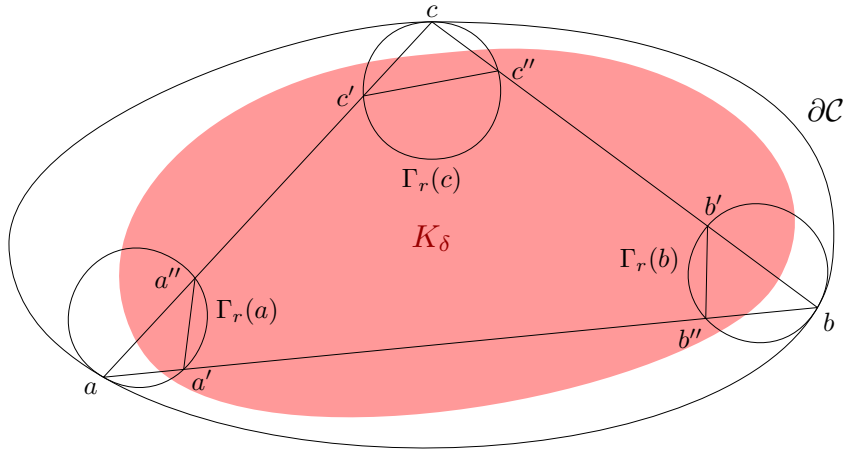}
  \caption{Premier cas \label{dessinA}}
\end{figure}

\begin{proof}[du théorème~\ref{upbound}] 
  Soit $T = a b c$ un triangle idéal de $(\mathcal{C} ,
  d_{\mathcal{C}})$ dont on oriente les sommets dans le sens
  trigonométrique.
  
  On va traiter trois cas selon que la longueur euclidienne des côtés
  du triangle $T$ est ou non inférieure à $r$.
  
\smallskip
\noindent  \textbf{Cas 1 :}~$d(a , b) \geq r$, $d(b , c) \geq r$ et $d(a , c) \geq r$
  (figure~\ref{dessinA}).
  
  En considérant le sommet $a$, désignons par $a'$ et $a''$ les points
  d'intersection du cercle $\Gamma_{r}(a)$ avec les segments $]a , b[$ et
  $]a , c[$ respectivement.  Ainsi le triangle $a a' a''$ a ses
  sommets orientés dans le sens trigonométrique et a $\Gamma_{r}(a)$ pour
  cercle euclidien circonscrit.  En procédant de même avec les sommets
  $b$ et $c$, on obtient les triangles $b b' b''$ et $c c' c''$.
  
  D'après le lemme~\ref{Distbord}~(ii), on a ici $d(a' , \mathcal \partial
  C) \geq \delta$ et $d(a'' , \mathcal \partial C) \geq \delta$, d'où $a', a'' \in
  K_{\delta}$ et par suite $[a' , a''] \subset K_{\delta}$ puisque $K_{\delta}$ est
  convexe.  De même, on a $[b' , b''] \subset K_{\delta}$ et $[c' , c''] \subset
  K_{\delta}$.
  
  Cela entraîne que $T \cap V_{\delta}$ est contenu dans la réunion des
  triangles $a a' a''$, $b b' b''$ et $c c' c''$, chacun d'eux étant
  d'aire majoré par $\pi$ en comparaison avec la géométrie hyperbolique
  associée aux disques ouverts $D_{r}(a)$, $D_{r}(b)$ et $D_{r}(c)$.
  On en déduit donc que $\mu_{\mathcal C}(T) \leq 3 \pi + \mu_{\mathcal
    C}(K_{\delta})$.
  
\smallskip  
\noindent  \textbf{Cas 2 :}~$d(a , b) \leq r$ et $d(b , c) \leq r$.
  
  Soit $m$ le projeté orthogonal du point $b$ sur la droite $(a c)$.
  Rappelons que $S(a,b)$ désigne le rectangle de base $[a,b]$ et
de hauteur $r$ donné par le lemme \ref{Recouvrement}.
 
  \begin{itemize}
  \item[\textbullet] Supposons que $m$ appartienne au segment $[a , c]$
  (figure~\ref{dessinB}).

\begin{figure}[hbtp] 
  \centering \includegraphics{dessinA-2}
\caption{Cas $m \in [a , c]$\label{dessinB}}
\end{figure}

Puisque le triangle $a b m$ est rectangle en $m$, on a $d(m , b) \leq
d(a , b)$.  Si $p$ est le projeté orthogonal de $m$ sur la droite $(a
, b)$, on a de même $d(m , p) \leq d(m , b)$ et par suite $d(m , p) \leq
d(a , b) \leq r$.  Comme en outre $p \in [a , b]$, il en résulte que $m$
appartient au rectangle fermé $\overline{S(a , b)}$.  Le même
raisonnement avec le triangle $b c m$ montre que $m$ appartient
également au rectangle fermé $\overline{S(b , c)}$, ce qui
entra\^{\i}ne que les triangles $a b m$ et $b c m$ sont inclus
respectivement dans les convexes $S(a , b)$ et $S(b , c)$.  Or
$\overline{a b m} \cup \overline{b c m} = \overline{a b c}$ puisque $m
\in [a , c]$, d'où $\overline{T} = \overline{a b c} \subset \overline{S(a ,
  b)} \cup \overline{S(b , c)}$.  On en déduit ainsi que $\mu_{\mathcal
  C}(T) \leq \mu_{\mathcal C}(S(a , b)) + \mu_{\mathcal C}(S(b , c))$ et
par conséquent $\mu_{\mathcal C}(T) \leq 4 \pi E\Bigl( \frac{2 R}{r}\Bigr)$
d'après le lemme~\ref{Recouvrement}.

\item[\textbullet] Supposons que $a$ appartienne au segment $[m , c]$
(figure~\ref{dessinC}).

\begin{figure}[hbtp]
  \centering \includegraphics{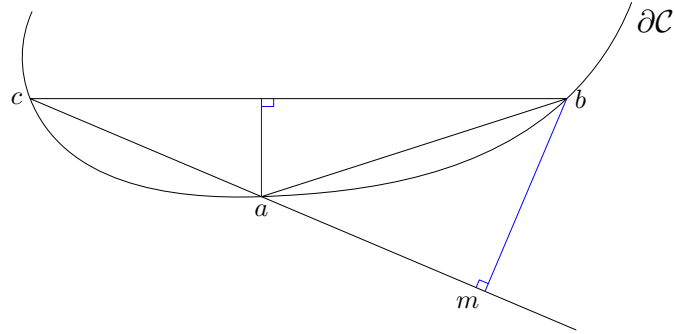}
\caption{Cas $a \in [m , c]$\label{dessinC}}
\end{figure}

Alors $d(a , c) \leq d(m , c)$, ce qui avec $d(m , c) \leq d(b , c)$ (le
triangle $b c m$ étant rectangle en $m$) conduit à $d(a , c) \leq d(b ,
c) \leq r$.  En outre, on a ici que, en notant $\hat{a}$ l'angle au sommet $a$ du triangle $abc$,
$\hat{a} \geq \pi / 2$ et par suite le
projeté orthogonal de $a$ sur la droite $(b c)$ est dans le segment
$[b , c]$, ce qui permet de conclure comme dans le point précédent.

\item[\textbullet] Supposons que $c$ appartienne au segment $[a , m]$.

On applique alors le point précédent en échangeant les rôles de $a$ et
$c$.

\end{itemize}

\smallskip
\noindent\textbf{Cas 3:}~$d(a , b) \leq r$, $d(b , c) \geq r$ et $d(a , c) \geq r$.

C'est la situation la plus délicate à traiter.  Tout comme dans le
premier cas, désignons par $c'$ et $c''$ les points d'intersection du
cercle $\Gamma_{r}(c)$ avec les segments $]a , c[$ et $]b , c[$
respectivement.

Dans ce qui suit, $\hat(a)$ et $\hat(b)$ sont les angles en $a$ et $b$
du triangle $abc$.

\begin{itemize}
\item[\textbullet] Supposons $\hat{a} \leq \pi / 2$ et $\hat{b} \leq \pi / 2$
(figure~\ref{dessinD}).

Soient $p$ et $q$ les sommets du rectangle $S(a , b)$ autres que $a$
et $b$ tels que $p - q = a - b$.  D'après le lemme~\ref{Rectangles},
la distance euclidienne de $p$ au centre du cercle $\Gamma_{r}(a)$ est
inférieure ou égale à $(3 / 4) r$, ce qui montre que $d(p , \Gamma_{r}(a))
\geq r / 4$, d'où $d(p , \partial \mathcal C) \geq r / 4$ puisque $p \in D_{r}(a)
\subset \mathcal C$ (et donc $d(p , \partial \mathcal C) \geq d(p , \Gamma_{r}(a))$).
De même, en considérant $\Gamma_{r}(p)$, on a $d(q , \partial \mathcal C) \geq r /
4$.  De $r / 4 \geq \delta$, on déduit alors que $p$ et $q$ sont dans le
convexe $K_{\delta}$ et par suite $[p , q] \subset K_{\delta}$.

Par ailleurs, comme $\hat{a} \leq \pi / 2$ (resp. $\hat{b} \leq \pi / 2$),
la droite $(a c)$ (resp. $(b c)$) coupe le segment $[p , q]$
(parallèle à $(a b)$ qui n'est parallèle ni à $(a c)$, ni à $(b c)$)
en un unique point $m_{a}$ (resp. $m_{b}$).  On a donc l'adhérence de
$T = a b c$ qui est incluse dans la réunion de l'adhérence du triangle
$c c' c''$ et des enveloppes convexes de $\{ a , b , m_{a} , m_{b} \}$
et $\{ m_{a} , m_{b} , c' , c'' \}$.

\begin{figure}[hbtp]
  \centering \includegraphics{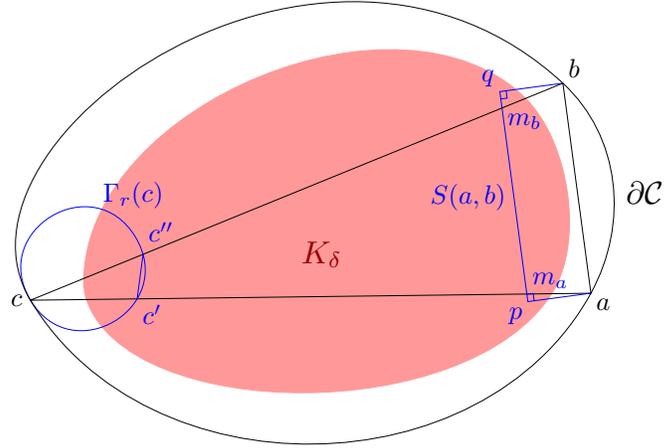}
\caption{Cas $\hat a \leq \pi / 2$ et $\hat b \leq \pi / 2$ \label{dessinD}}
\end{figure}

Comme d'une part $m_{a}, m_{b} \in K_{\delta}$ et $c', c'' \in K_{\delta}$ (même
raison que dans le premier cas), l'enveloppe convexe de $\{ m_{a} ,
m_{b} , c' , c'' \}$ est dans $K_{\delta}$ et puisque d'autre part $a, b,
m_{a}, m_{b}$ appartiennent au rectangle fermé $\overline{S(a , b)}$
qui est convexe, l'enveloppe convexe de $\{ a , b , m_{a} , m_{b} \}$
est dans $\overline{S(a , b)}$.

Il en résulte que
$$
\mu_{\mathcal C}(T) \leq \mu_{\mathcal C}\Bigl(\!S(a , b)\!\Bigr) +
\mu_{\mathcal C}(K_{\delta}) + \mu_{\mathcal C}(c c' c'')
$$
et par conséquent
$$
\mu_{\mathcal C}(T) \leq \pi \biggl( \!2 E\Bigl(\frac{2 R}{r}\Bigr) + 1 \!\biggr) + \mu_{\mathcal C}(K_{\delta})
$$
en vertu du lemme~\ref{Recouvrement}.

\item[\textbullet] Supposons $\hat{a} > \pi / 2$ (le cas $\hat{b} > \pi / 2$ se traite
de fa\c{c}on similaire) (figure~\ref{dessinE}).

Introduisons comme précédemment les points $p$, $q$ et $m_{b}$
(puisque $\hat{b} \leq \pi / 2$) et soit en outre $a''$ le point
d'intersection du cercle $\Gamma_{r}(a)$ avec le segment $]a , c[$.

\begin{figure}[hbtp]
  \centering \includegraphics{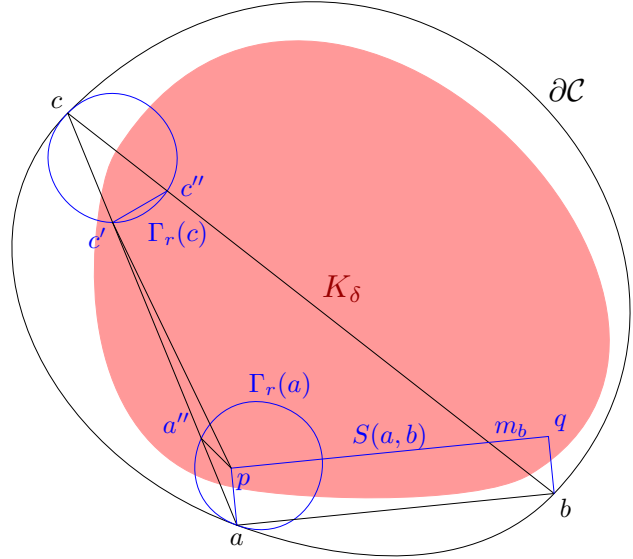}
  \caption{Cas où $\hat a > \pi / 2$ \label{dessinE}}
\end{figure}

Sachant que $d(a , c) \geq r$, on a $d(a'' , \mathcal \partial C) \geq \delta$
d'après le lemme~\ref{Cordes}~(ii), d'où $a'' \in K_{\delta}$.  D'autre
part, en raisonnant comme au point précédent, on a $p, m_{b}, c', c''
\in K_{\delta}$ et $a, b, p, m_{b} \in \overline{S(a , b)}$.  Par suite,
l'enveloppe convexe de $\{ p , m_{b} , c' , c'' \}$ et le triangle $p
a'' c'$ sont inclus dans $K_{\delta}$ alors que l'enveloppe convexe de $\{
a , b , p , m_{b} \}$ est dans le rectangle fermé $\overline{S(a ,
  b)}$.

Enfin, la convexité de $D_{r}(a)$ et le fait que $a, a'' , p \in
\overline{D_{r}(a)}$ assurent que le triangle $a a" p$ est inclus dans
$D_{r}(a)$.

Puisque l'adhérence de $T = a b c$ est contenue dans la réunion des
adhérences des triangles $p a'' c'$, $a a" p$ et $c c' c''$ ainsi que
des enveloppes convexes de $\{ p , m_{b} , c' , c'' \}$ et $\{ a , b ,
p , m_{b} \}$, il s'ensuit que
$$
\mu_{\mathcal C}(T) \leq \mu_{\mathcal C}\bigl(S(a , b)\bigr) +
\mu_{\mathcal C}(K_{\delta}) + \mu_{\mathcal C}(a a'' p) + \mu_{\mathcal C}(c
c' c'')\text{,}
$$
d'où
$$
\mu_{\mathcal C}(T) \leq 2 \pi \biggl( \!E\Bigl(\frac{2 R}{r}\Bigr) +1 \!\biggr) + \mu_{\mathcal C}(K_{\delta})
$$
d'après le lemme~\ref{Recouvrement}.
\end{itemize}
\end{proof}

\begin{rem}
  Questions ouvertes concernant les majorants.
 
  Pour conclure le présent travail, tentons de donner les hypothèses
  les plus faibles que l'on doit imposer à un domaine convexe donné
  $\mathcal{C}$ du plan pour espérer obtenir une majoration de l'aire
  des triangles idéaux de $(\mathcal C , d_{\mathcal C})$.
  
  \makeatletter \renewcommand\theenumi{\@arabic\c@enumi} \makeatother
\begin{enumerate}
\item D'après le corollaire~\ref{plat}, il ne peut y avoir de segment
  ouvert dans le bord $\partial \mathcal{C}$, ce qui se traduit par le fait
  que $\mathcal{C}$ doit être \textit{affinement strictement convexe},
  autrement dit que tout segment de droite ouvert entre deux points de
  $\partial \mathcal{C}$ est contenu dans $\mathcal C$.
  
\item D'après le corollaire~\ref{coin}, il ne peut y avoir de
  \og\,coin\, \fg{} dans $\mathcal C$, c'est-à-dire de point de $\partial
  \mathcal{C}$ en lequel $\mathcal{C}$ admet deux droites d'appui
  distinctes.  Du point de vu analytique, ceci signifie que le bord
  $\partial \mathcal{C}$ doit être localement le graphe d'une fonction
  convexe partout dérivable. Mais avec l'hypothèse de convexité, cela
  implique que le bord de $\mathcal{C}$ est une courbe de classe $C^1$
  (voir \cite{bourbaki}, I.32, § 4).  Ainsi, on peut se concentrer sur
  un domaine $\mathcal{C}$ affinement strictement convexe dont le bord
  est $C^1$.
     
\item À la proposition~\ref{prop:airefinie}, nous avons cependant eu
  besoin d'avoir $\partial \mathcal{C}$ de classe $C^{2}$ pour montrer que
  l'aire des triangles idéaux de $(\mathcal C , d_{\mathcal C})$ est
  finie, et d'ajouter l'hypothèse que la courbure de $\partial \mathcal{C}$
  n'est jamais nulle pour exhiber un majorant de cette aire.
   \end{enumerate}
   \makeatletter \renewcommand\theenumi{\@roman\c@enumi} \makeatother
   
   Ainsi, lorsque $\mathcal{C}$ est affinement strictement convexe
   avec un bord de classe $C^{1}$ sans être $C^{2}$, la finitude de
   l'aire des triangles idéaux de $(\mathcal C , d_{\mathcal C})$
   reste un problème ouvert.  Tout comme l'est la question de savoir
   s'il existe un majorant de l'aire de ces triangles lorsque
   $\mathcal{C}$ est affinement strictement convexe avec un bord de
   classe $C^{2}$ dont la courbure s'annule en certains points.
 
\end{rem}

  \par
  \setcounter{annexe}{0}%
  \setcounter{subsection}{0}%
  \makeatletter \renewcommand\theannexe{\@Alph\c@annexe}
  \renewcommand\thesubsection{\@Alph\c@annexe.\@arabic\c@subsection}
  \renewcommand\thelemme{\@Alph\c@annexe.\@arabic\c@lemme}
  \makeatother

\annexe{Les domaines triangulaires} \label{annexeA}

Rappelons que l'on s'est donné un domaine triangulaire $\Delta=mpq \subset \R^2$
et un triangle idéal $T=abc$ de $(\Delta , d_{\Delta})$.
Il existe donc $\lambda, \mu, \nu \in \, ]0 , 1[$ tels que
\begin{align*}
  a&= (1 - \lambda) m + \lambda p,& b &= (1 - \mu) p + \mu q & &\text{et}&  c &=
  (1 - \nu) q + \nu m \text{.}
\end{align*}

On considère par ailleurs le domaine triangulaire $\Delta_{0}$ de $\R^2 \subset \R^3$
dont les sommets
sont $0 = (0 , 0 , 0)$, $e_{1} = (1 , 0 , 0)$ et $e_{2} = (0, 1 , 0)$ et le
triangle idéal $T(\alpha)$ de $(\Delta_{0} , d_{\Delta_{0}})$ ayant pour sommets
$a(\alpha) = (\alpha , 1 - \alpha , 0)$, $b(\alpha) = (0 , 1 - \alpha , 0)$ et $c(\alpha) = (\alpha , 0 , 0)$.

\subsection{Preuve du lemme $\ref{etape1}$} \label{annexeA1}
Introduisons d'abord  $f \: \R^{2} \to
\R^{3}$ l'application affine injective qui envoie les points
$m$, $p$ et $q$ sur $e_{1}$, $e_{2}$ et $e_{3} = (0 , 0 , 1)$ respectivement.
En notant
\begin{align*}
  a'& = (1 - \lambda) e_{1} + \lambda e_{2}, &b'& = (1 - \mu) e_{2} + \mu e_{3}
  &&\text{et} &c'& = (1 - \nu) e_{3} + \nu e_{1},
\end{align*}
le triangle $T' = a' b' c'$ est alors l'image de $T$ par $f$ et a ses
sommets sur les côtés du triangle $\Delta' = f(\Delta)$.

Considérons ensuite les points
\begin{eqnarray*}
a'' & = & \alpha (u e_{1}) + (1 - \alpha) (v e_{2}), \\
b'' & = & (1 - \alpha) (v e_{2}) + \alpha (w e_{3}), \\
c'' & = & (1 - \alpha) (w e_{3}) + \alpha (u e_{1})
\end{eqnarray*}
avec
\begin{eqnarray*}
\alpha & = & \frac{\lambda \mu \nu}{(1-\lambda)(1-\mu)(1-\nu)+\lambda\mu\nu}
\in \, ]0 , 1[, \\
u & = & \frac{(1 - \lambda) (1 - \mu)}{\lambda \mu}w \neq 0
\quad \textrm{et} \quad
v = \frac{\lambda \nu}{(1 - \lambda) (1 - \nu)}w \neq 0.
\end{eqnarray*}
Remarquons que, quitte à remplacer $\lambda$ par $1 - \lambda$,
$\mu$ par $1 - \mu$ et $\nu$ par $1 - \nu$, on peut supposer $\alpha\in\,]0,1/2]$.
On vérifie sans peine que $a''$ est sur la droite $(0 a')$, $b''$ sur
la droite $(0 b')$ et $c''$ sur la droite $(0 c')$, de sorte qu'en
désignant par $\pi : \R^{3} \to \mathbf{P}^{2}(\R)$ la
projection canonique, on obtient
\begin{equation} \label{eq:tideal}
\pi(T') = \pi(T'')
\end{equation}
ainsi que
\begin{equation} \label{eq:triangle}
\pi(\Delta') = \pi(\Delta''),
\end{equation}
où $T'' = a'' b'' c''$ et $\Delta''$ est le triangle de sommets $u e_{1}$,
$v e_{2}$ et $w e_{3}$.

Enfin, si $L$ est l'application linéaire surjective de
$\R^{3}$ sur $\R^{2}$ qui envoie respectivement $u
e_{1}$, $v e_{2}$ et $we_{3}$ sur $e_{1}$, $e_{2}$ et $0$, alors
$L(\Delta'') = \Delta_{0}$ et $L(T'') = T(\alpha)$.

Comme les applications $f$ et $L$ sont affines et qu'on a les
égalités~(\ref{eq:tideal}) et (\ref{eq:triangle}), il en résulte que
les géométries de Hilbert $(\Delta , d_{\Delta})$ et $(\Delta_{0} , d_{\Delta_{0}})$
sont isométriques avec correspondance entre les triangles idéaux $T$
et $T(\alpha)$.

\begin{figure}[htpb]
  \centering \includegraphics{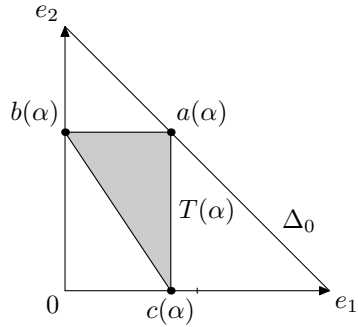}
  \caption{Triangle idéal $T(\alpha)$ pour le domaine triangulaire $\Delta_0$}
\end{figure}

\subsection{L'aire des triangles idéaux pour un domaine triangulaire}\label{annexeA2}
L'aire du triangle idéal $T(\alpha)$ vaut

\begin{equation*} \label{eq00}
  \mathcal{A}(\alpha) = \frac{\pi}{12}\int_0^\alpha \!\int_{(1-\alpha)(1-x/\alpha)}^{1-\alpha}
  \frac{\ed{x}\ed{y}}{xy(1-x-y)}
\end{equation*}

et se décompose comme suit~:

\begin{eqnarray}
  \frac{12}{\pi}\mathcal{A}(\alpha)&=&\int_0^\alpha \!\int_{(1-\alpha)(1-x/\alpha)}^{1-\alpha}
  \Bigl(\frac{1}{y}+\frac{1}{1-x-y}\Bigr)\frac{\ed{x}\ed{y}}{x(1-x)} \nonumber \\
  &=&-2 \!\int_0^\alpha \frac{\ln(1-x/\alpha)}{x(1-x)}\ed{x}  \label{eq01} \\
  &&+\int_0^\alpha \ln\left( \Bigl( \!\frac{1-2\alpha}{\alpha^2} \!\Bigr) x +1 \right) \frac{\ed{x}}{x(1-x)}~. \label{eq02}
\end{eqnarray}

Le cacul de l'intégrale~(\ref{eq01}) donne

\begin{eqnarray}
\int_0^\alpha \frac{\ln(1-x/\alpha)}{x(1-x)}\ed{x} &=&
\underbrace{\int_0^\alpha \frac{\ln(1-x/\alpha)}{x} \ed{x}}_{\text{On pose}\ u~=~x/\alpha} %
 +\underbrace{\int_0^\alpha \frac{\ln(1-x/\alpha)}{1-x} \ed{x}}_{\text{Intégration par parties}} \nonumber \\
&=& \int_0^1 \frac{\ln(1-u)}{u} \ed{u}
- \underbrace{\int_0^\alpha  \frac{\ln\Bigl(1+\frac{\alpha-x}{1-\alpha}\Bigr)}{\alpha-x} \ed{x}}
_{\text{On pose}\ v~=~\frac{\alpha-x}{1-\alpha}} \nonumber \\
&=&\int_0^1 \frac{\ln(1-u)}{u} \ed{u} - \int_0^{\frac{\alpha}{1-\alpha}} \frac{\ln(1+v)}{v} \ed{v}. \label{eq03}
\end{eqnarray}

Par ailleurs, l'intégrale~(\ref{eq02}) s'écrit

\begin{eqnarray}
  \int_0^\alpha\frac{\ln\Bigl(\frac{1-2\alpha}{\alpha^2} x +1 \Bigr)}{x(1-x)}\ed{x} & = & %
  \underbrace{\int_0^\alpha\frac{\ln\Bigl(\frac{1-2\alpha}{\alpha^2}x +1\Bigr)}{x}\ed{x}}
  _{\text{On pose } u~=~\left( \frac{1-2\alpha}{\alpha^2} \right)x}
  + \underbrace{\int_0^\alpha\frac{\ln\Bigl(\frac{1-2\alpha}{\alpha^2}x +1 \Bigr)}{1-x}\ed{x}}
  _{\text{Intégration par parties}} \nonumber \\
  & = & \int_0^{\frac{1-2\alpha}{\alpha}} \frac{\ln(1+u)}{u} \ed{u} %
  + \underbrace{\int_0^\alpha\frac{\ln\Bigl(1+\frac{\alpha-x}{1-\alpha}\Bigr)}
  {1+\frac{1-2\alpha}{\alpha^2}x}\ed{x}}_{\text{On pose } v~=~\frac{\alpha-x}{1-\alpha}} \nonumber \\
  & = & \int_0^{\frac{1-2\alpha}{\alpha}} \frac{\ln(1+u)}{u} \ed{u} + {\scriptstyle \frac{1-2\alpha}{\alpha}} 
\!\! \int_0^{\frac{\alpha}{1-\alpha}} \frac{\ln(1+v)}{1-\frac{1-2\alpha}{\alpha}v} \ed{v}.
  \label{eq04}
\end{eqnarray}

On utilise alors dans (\ref{eq03}) et (\ref{eq04}) la nouvelle variable $t = \frac{1-2\alpha}{\alpha}$ qui parcourt $[0 , +\infty[$
lorsque $\alpha$ décrit $]0 , 1/2]$.  En posant $\mathcal{F}(t)= \frac{12}{\pi}
\mathcal{A}(\alpha)$, on obtient ainsi

\begin{multline} \label{eq05}
  \mathcal{F}(t)= -2 \!\! \int_0^1\frac{\ln(1-u)}{u} \ed{u} + 2 \!\! \int_0^{\frac{1}{1+t}} \frac{\ln(1+v)}{v} \ed{v} \\
  + \int_0^t \frac{\ln(1+u)}{u} \ed{u} + t \!\! \int_0^{\frac{1}{1+t}}
  \frac{\ln(1+v)}{1-tv} \ed{v}.
\end{multline}

Enfin, le changement de variable $w = \frac{1-tv}{1+t}$ dans la dernière intégrale
de (\ref{eq05}) conduit à

\begin{multline*}
  \mathcal{F}(t) = -2 \!\! \int_0^1 \frac{\ln(1-u)}{u} \ed{u} + 2 \!\! \int_0^{\frac{1}{1+t}} \frac{\ln(1+v)}{v} \ed{v} \\
  + \int_0^t \frac{\ln(1+u)}{u} \ed{u}
  + \int_{\frac{1}{(1+t)^2}}^{\frac{1}{1+t}} \! \frac{\ln(1-w)}{w} \ed{w} %
  + \ln\!\left( \frac{1+t}{t} \right)\!\ln(1+t).
\end{multline*}

La fonction $\mathcal{F}$ se dérive sans trop de difficultés et, après
simplifications, on obtient finalement

\begin{equation*}
  \mathcal{F}'(t) = \frac{4}{1+t}\ln\!\left( \frac{(1+t)^2}{(1+t)^2-1} \right) > 0.
\end{equation*}

On en déduit donc que $\mathcal{F}$ est strictement croissante sur $[0, +\infty[$
et par suite son minimum est atteint en $t = 0$ seulement.
Autrement dit, l'application $\mathcal{A} : \, ]0 , 1/2] \longrightarrow
\R$ est strictement décroissante et atteint son minimum en le
seul point $\alpha = 1/2$.

\annexe{Lemmes techniques du théorème~\ref{upbound}} \label{annexeB}

\begin{lemme} \label{Cordes} 
  Étant donné des réels $0 < \varrho < \varrho'$, soient $\Gamma_{\varrho}$ et $\Gamma_{\varrho'}$
  les cercles euclidiens de $\R^{2}$ passant par l'origine et
  de centres respectifs $c = (0 , \varrho)$ et $c' = (0 , \varrho')$.  Pour
  chaque point $m$ du segment $[0 , c]$ et pour chaque vecteur non nul
  $v \in \R^{2}$, on note $p$ (resp. $p'$) le point
  d'intersection de la demi-droite fermée $m + \R_{-} v$ avec
  $\Gamma_{\varrho}$ (resp. $\Gamma_{\varrho'}$) et $q$ (resp. $q'$) le point
  d'intersection de la demi-droite fermée $m + \R_{+} v$ avec
  $\Gamma_{\varrho}$ (resp.  $\Gamma_{\varrho'}$).

  On a alors
  $$
  d(p , q) \geq \left( \frac{\varrho}{\varrho'} \right) \!\! d(p' , q').
  $$
  En outre, lorsque $m = 0$ (et donc $p = p' = 0$), on a
  $$
  d(q , \Gamma_{\varrho'}) \geq \left( \frac{\varrho' - \varrho}{2 \varrho \varrho'} \right) \!\! d(0 , q)^{2}.
  $$
\end{lemme}

\begin{figure}[hbtp]
  \centering \includegraphics{dessin1-1}
  \caption{Lemme~\ref{Cordes}\label{dessin1}}
\end{figure}

\begin{proof}[du lemme~\ref{Cordes}]

  Lorsque $p = q$, la droite $m + \R v$ est égale à
  $\R × \{ 0 \}$, ce qui entraîne que $p' = q'$ et par suite
  le lemme est trivialement vérifié.

  Supposons donc $p$ et $q$ distincts. Soient $a$ et $a'$ les milieux de $(p , q)$ et $(p' , q')$
  respectivement.  Comme $c$ (resp. $c'$) est sur la médiatrice de $(p
  , q)$ (resp. $(p' , q')$), le vecteur $a - c$ (resp. $a' - c'$) est
  orthogonal à $v$.  Il en résulte que $a - c$ et $a' - c'$ sont
  colinéaires, d'où l'existence d'un réel  $\lambda$ tel que $a - c = \lambda
  (a' - c')$.

  Or $m$, $c$ et $c'$ étant alignés ainsi que $m$, $a$ et $a'$, on a
  aussi (Thalès) $m-c=\lambda(m-c')$.  En écrivant $m=tc$ avec $t \in [0 ,
  1]$ et sachant que $c'= (\varrho'/\varrho)c$, il vient
  $$\lambda = \frac{1 - t}{(\varrho' / \varrho) - t} \in [0 , 1[$$
  et par conséquent
  $(\varrho' / \varrho) \lambda \in [0 , 1]$.  On en déduit que
\begin{eqnarray*}
  d(a , q)^{2} = {\varrho}^{2} - d(a , c)^{2}
 & =& {\varrho}^{2} - \lambda^{2} d\bigl(a' , c'\bigr)^{2}\\
  & = & \Bigr(\frac{\varrho}{\varrho'}\Bigl)^{\! 2} \! \biggl\{ {\varrho'}^{2}
  - \Bigl(\frac{\varrho'}{\varrho} \lambda\Bigr)^{\! 2} d\bigl(a' , c'\bigr)^{2} \biggr\} \\
  & \geq & \Bigr(\frac{\varrho}{\varrho'}\Bigl)^{\! 2} \! \bigl({\varrho'}^{2} -
  d\bigl(a' , c'\bigr)^{2}\bigr)
  = \Bigr(\frac{\varrho}{\varrho'}\Bigl)^{\! 2} \! d\bigl(a' , q'\bigr)^{2}\text{,}
\end{eqnarray*}
ce qui démontre la première inégalité.

En ce qui concerne la deuxième inégalité, on a
$d\bigl(q ,\Gamma_{\varrho'}\bigr) = \varrho' - d\bigl(q , c'\bigr)$ et
$$
{\varrho'}^{2} - d\bigl(q , c'\bigr)^{2} = 2 \varrho' d(0 , q) \cos{\hat{0}}
- d(0 , q)^{2}
$$
dans le triangle $0 q c'$,
ce qui entraîne que
$$
d\bigl(q , \Gamma_{\varrho'}\bigr) \Bigl(d\bigl(q , c'\bigr) + {\varrho'}\Bigr) =
2 \varrho' d(0 , q) \cos{\hat{0}} - d(0 , q)^{2}\text{.}
$$
Par ailleurs, dans le triangle isocèle $0 q c$, on a $d(0 , q) = 2
\varrho \cos{\hat{0}}$, d'où il résulte que
$$
d\bigl(q , \Gamma_{\varrho'}\bigr) \Bigl(d\bigl(q , c'\bigr) + {\varrho'}\Bigr) =
\bigl((\varrho' / \varrho) - 1\bigr) d(0 , q)^{2}\text{.}
$$
En remarquant alors que $d\bigl(q , c'\bigr) \leq \varrho'$, on obtient la
relation désirée.  \end{proof}

\begin{lemme} \label{Rectangles} 
  Soient $r > 0$ fixé, $\Gamma_{r}$ le cercle euclidien de
  $\R^{2}$ passant par l'origine et de centre $c = (0 , r)$ et
  $D_{r}$ le disque ouvert correspondant.  Pour tout $h \in [0 , r]$,
  notons $p = (-\alpha , h)$ (resp. $q = (\alpha , h)$) l'intersection de la
  droite d'équation $y = h$ avec $\displaystyle \Gamma_{r} \cap
  (\R_{-} × \R)$ (resp. $\displaystyle \Gamma_{r} \cap
  (\R_{+} × \R)$) et soit $p' = (-\alpha , h + r)$ (resp.
  $q' = (\alpha , h + r)$).

  Alors, si $d(p , q) = 2 \alpha \leq r$, on a $d(c , p') = d(c , q') \leq (3
  / 4) r$ (les points $p'$ et $q'$ sont donc en particulier dans
  $D_{r}$).
\end{lemme}

\begin{figure}[hbtp]
  \centering \includegraphics{dessin1-2}
  \caption{Lemme~\ref{Rectangles} \label{dessin2}}
\end{figure}

\begin{proof}[du lemme~\ref{Rectangles}]
  \begin{align}
  &\text{On a}& d\bigl(c , q'\bigr)^{2} &=  \alpha^{2} + h^{2} \nonumber& \\
&\text{et} & r^{2} &=  d(c , q)^{2} = \alpha^{2} + (r - h)^{2}\text{,} &\label{eq:rec}
\end{align}
d'où il résulte que $d\bigl(c , q'\bigr)^{2} = 2 r h$. Comme $\alpha \leq r
/ 2$, on déduit de (\ref{eq:rec}) que $h \leq \bigl( 1 - \sqrt{3} / 2
\bigr) r$ et par suite
$$
d\bigl(c , q'\bigr)^{2} \leq \Bigl( 2 - \sqrt{3} \Bigr) r^{2} \leq
r^{2} / 2\text{.}
$$
Donc $d\bigl(c , q'\bigr) \leq \bigl( \sqrt{2} / 2 \bigr) r \leq \bigl(3 /4\bigr) r$.
\end{proof}

\begin{proof}[du lemme~\ref{Distbord}]\null

  Il existe au moins un point d'intersection $a'$ entre $]a , b[$ et
  $\Gamma_{r}(a)$ car sinon soit $b$ est dans $\Gamma_{r}(a)$, ce qui
  contredit le fait que le cercle de rayon $2 r$ roule à l'intérieur
  de $\overline{\mathcal C}$, soit $b$ est dans le demi-plan fermé
  bordé par la tangente à $\partial \mathcal C$ en $a$ qui ne contient pas
  le convexe strict $\mathcal C$, ce qui est là encore impossible.
  L'unicité résulte du fait qu'un cercle coupe une droite en au plus
  deux points distincts et il y a déjà $a$ et $a'$ dans l'intersection
  de $\Gamma_{r}(a)$ avec la droite $(a b)$.

  Puisque $a' \in D_{2 r}(a) \subset \mathcal{C}$, on a
  $$
  d\bigl(a' , \partial \mathcal C\bigr) \geq d\bigl(a' , \Gamma_{2
    r}(a)\bigr)\text{.}
  $$
  Or, d'après la deuxième inégalité du lemme~\ref{Cordes} avec $\varrho
  = r$ et $\varrho' = 2 r$, on a
  $$
  d\bigl(a' , \Gamma_{2 r}(a)\bigr) \geq \frac{1}{4 r} d\bigl(a ,
  a'\bigr)^{2}
  $$
  et d'après la première inégalité de ce même lemme avec $\varrho = r$,
  $\varrho' = R$, on a $d(a , a') \geq (r / R) d(a , m)$, où $m$ est le point
  d'intersection autre que $a$ entre $\Gamma_{R}(a)$ et la droite $(a b)$.
  Comme $a \in \overline{D_{r}(a)} \subset
  \overline{\mathcal C} \subset \overline{D_{R}(a)}$,
  on a $b \in [a' , m] \subset [a , m]$, d'où $d(a , m) \geq d(a , b)$ et par suite
  $$
  d\bigl(a' , \partial \mathcal C\bigr) \geq \frac{1}{4 r}×
  \Bigl(\frac{r}{R}\Bigr)^{2} d(a , b)^{2} = \frac{r}{4 R^{2}} d(a ,
  b)^{2}\text{.}
  $$
  \end{proof}

En ce qui concerne le lemme~\ref{Recouvrement}, il va se déduire du lemme
technique suivant~:

\begin{lemme}\label{recouvrebis}
Étant donné un convexe strict $\mathcal{C}$ du plan, soient $a$ et $b$ deux points distincts
de $\partial \mathcal C$ tels que $d(a , b) \leq r$. Alors :
\begin{enumerate}
\item Pour chaque $m \in \, ]a , b[$, il existe
$\omega_m \in \partial {\mathcal C} \!\! \setminus \!\! \{ a , b \}$ tel que $d(m , \partial {\mathcal C}) = d(m , \omega_m)$
avec $m - \omega_m \perp \partial {\mathcal C}$.

\item L'intersection $C(a , b)$, entre  $\partial \mathcal C$ et l'un des deux demi plan 
fermés $H^{-}(a , b)$ et $H^{+}(a , b)$ de $\R^{2}$ bordés par la droite $(ab)$, vérifie
$d(m , C(a , b)) = d(m , \partial {\mathcal C})$ quel que soit $m \in \, ]a , b[$.
\end{enumerate}
\end{lemme}

\begin{proof}[du lemme~\ref{recouvrebis}]
\begin{enumerate}
\item Pour $m \in \, ]a , b[$ fixé, la fonction $f : \partial {\mathcal C} \longrightarrow \R$ définie par $f(\omega) = d^2(m , \omega)$
étant continue sur le compact $\partial {\mathcal C}$, l'existence
de $\omega_m \in \partial {\mathcal C}$ tel que $d(m , \partial {\mathcal C}) = d(m , \omega_m)$ en découle.
De plus, comme $\partial {\mathcal C}$ et $f$ sont différentiables, $\omega_m$ est un point critique de $f$,
ce qui conduit à $m - \omega_m \perp \partial {\mathcal C}$.
Enfin, si on avait $\omega_m = a$, on aurait alors $(a b) \perp \partial {\mathcal C}$ en $a$ et par suite
$b \in D_r(a)$ puisque $d(a , b) \leq r$. Or $D_r(a) \subset {\mathcal C}$, d'où il résulterait que $b \in {\mathcal C}$,
ce qui est faux.
Par conséquent, on a $\omega_m \neq a$ ainsi que $\omega_m \neq b$ pour la même raison.

\item À présent, en notant $C^-(a , b) = H^-(a , b) \cap \partial {\mathcal C}$ et $C^+(a , b) = H^+(a , b) \cap \partial {\mathcal C}$,
supposons qu'il existe $m^-, m^+ \in [a , b]$  tels que
$d(m^- , C^-(a , b)) > d(m^- , \partial {\mathcal C})$ et $d(m^+ , C^+(a , b)) > d(m^+ , \partial {\mathcal C})$.

Soient alors $\omega^-, \omega^+ \in \partial {\mathcal C} \! \setminus \! \{ a , b \}$ tels que
$d(m^- , \partial {\mathcal C}) = d(m^- , \omega^-)$ et $d(m^+ , \partial {\mathcal C}) = d(m^+ , \omega^+)$
--- d'où nécessairement $\omega^- \in C^+(a , b)$ et $\omega^+ \in C^-(a , b)$ ---
avec $m^- - \omega^- \perp \partial {\mathcal C}$ et $m^+ - \omega^+ \perp \partial {\mathcal C}$.

On a donc $r \geq d(a , b) \geq d(m^- , b) \geq d(m^- , C^-(a , b)) > d(m^- , \partial {\mathcal C}) = d(m^- , \omega^-)$,
d'où $d(m^- , \omega^-) < r$ ainsi que $d(m^+ , \omega^-) < r$ de manière analogue.
De là, il résulte alors que le centre $c^- \in (\omega^- m^-)$ du cercle $\Gamma_r(\omega^-)$ est dans $H^-(a , b)$
et que le centre $c^+ \in (\omega^+ m^+)$ du cercle $\Gamma_r(\omega^+)$ est dans $H^+(a , b)$.

\begin{figure}[hbtp]
  \centering \includegraphics{dessinlf}
  \caption{Lemme~\ref{recouvrebis}~(i) \label{convexe1}}
\end{figure}

En désignant enfin par $S^-$ (resp. $S^+$) l'unique diamètre (segment fermé)
de $\Gamma_r(\omega^-)$ (resp. $\Gamma_r(\omega^+)$) parallèle à la droite $(a b)$,
l'enveloppe convexe de $S^- \cup S^+$ (parallélogramme plein) a une intersection
$[p , q]$ avec $(a b)$ telle que $d(p , q) \geq 2r$.
Or, la convexité de ${\mathcal C}$ implique que cette enveloppe convexe est incluse
dans $\bar {\mathcal C}$, et par suite $[p , q] \subset [a , b]$.
On a donc $d(a , b) \geq d(p , q) \geq 2r$, ce qui impossible puisque $d(a , b) \leq r$
par hypothèse.
\end{enumerate}
\end{proof}

\begin{proof}[du lemme~\ref{Recouvrement}]
\begin{itemize}
\item[\textbullet]
Pour chaque $m \in \, ]a , b[$, soit $\omega_m \in \partial {\mathcal C} \! \setminus \! \{ a , b \}$ tel que
$d(m , C(a , b)) = d(m , \omega_m)$ donné par le lemme~\ref{recouvrebis}.
Ceci entraîne que $d(m , \omega_m) \leq r$ et par suite le cercle $\Gamma_r(\omega_m)$ coupe
la droite $(a b)$ en deux points $p_m \in [a , m[$ et $q_m \in ]m , b]$ qui vérifient
$d(p_m , q_m) \geq (r / R) d(a , b)$ en vertu de la première inégalité du lemme~\ref{Cordes}
avec $\varrho = r$, $\varrho' = R$ et $v = b - a$.

Par ailleurs, d'après le lemme~\ref{Distbord}, soient $a'$ et $b'$ les uniques points d'intersection
de $]a , b[$ avec les cercles $\Gamma_r(a)$ et $\Gamma_r(b)$ respectivement, pour
lesquels on a $d(a , a') \geq (r / R) d(a , b)$ et $d(b' , b) \geq (r / R) d(a , b)$, toujours
d'après le lemme~\ref{Cordes}.

Alors $[a , a'[$, $]b' , b]$ et la famille $(]p_m , q_m[)_{m \in \partial {\mathcal C} \setminus \{ a , b \}}$
forment un recouvrement ouvert de $[a , b]$, dont on peut donc extraire 
un sous-recouvrement fini $\mathcal I$
qui est minimal pour l'inclusion.
En outre, d'après ce qui précède, les éléments de $\mathcal I$ sont des segments
de longueurs supérieures ou égales à $(r / R) d(a , b)$.

À présent, montrons qu'aucun point de $[a , b]$ ne peut appartenir
à plus de deux éléments de $\mathcal I$, c'est-à-dire 
que si un point $x_0 \in [a , b]$ vérifie $x_0 \in I \cap J$ avec $I, J \in \mathcal I$,
alors pour tout $H \in \mathcal I \! \setminus \! \{ I , J \}$, on a $x_0 \not \in H$.

En effet, après identification de $[a , b]$ avec un segment de $\mathbb{R}$ et quitte
à échanger les rôles de $I$ et $J$,
on a  $\inf(I) < \inf(J) < \sup(I) < \sup(J)$ par minimalité de $\mathcal I$.
Aussi, supposons qu'il existe $H \in \mathcal I \! \setminus \! \{ I , J \}$ tel que $x_0 \in H$.

Si on avait $\inf(J) \leq \inf(H)$, alors on aurait $\sup(J) < \sup(H)$ car $H \not \subset J$
(minimalité de $\mathcal I$) et par suite $J \subset I \cup H$, ce qui est faux puisque
$\mathcal I$ est minimal.
C'est donc que $\inf(H) < \inf(J)$.
Or ceci impose que $\sup(H) < \sup(J)$ car $J \not \subset H$,
d'où il vient que $\inf(H) < \inf(I)$ (sinon $H \subset I \cup J$) et donc
que $\sup(I) < \sup(H)$ (sinon $H \subset I$).
Mais alors, c'est que $I \subset H$, ce qui, là encore, est impossible.

Notons alors $\mathcal I = \{ I_1 , \ldots , I_n \}$, $A_2=\bigcup_{1\leq i<j\leq n}I_i\cap I_j$ et $A_1=[a,b]\setminus A_2$.
Les ensembles $A_1$ et $A_2$
forment ainsi une partition mesurable de $[a , b]$ et vérifient
$\textrm{Long}(I_1) + \dots + \textrm{Long}(I_n) = \textrm{Long}(A_1) + 2 \textrm{Long}(A_2)$.

Or, $\textrm{Long}(A_1) + 2 \textrm{Long}(A_2)
= (\textrm{Long}(A_1) + \textrm{Long}(A_2)) + \textrm{Long}(A_2)
\leq \textrm{Long}([a , b]) + \textrm{Long}([a , b]) = 2 d(a , b)$,
d'où il résulte que
$n (r / R) d(a , b) \leq 2 d(a , b)$
puisqu'on a vu que
$\textrm{Long}(I_k) \geq (r / R) d(a , b)$ pour tout $k = 1, \ldots , n$.

Conclusion~: $n \leq E(2 R / r)$.

\item[\textbullet]
On peut maintenant terminer la preuve du lemme~\ref{Recouvrement}
en remarquant que le rectangle fermé
$\overline{S(a , b)}$ est la réunion des rectangles fermés
$S_k$ de base le segment $I_k$ et de hauteur $r$ pour $1 \leq k \leq n$,
chacun d'eux
étant inclus dans un disque $\overline{D_{r}(\omega_{k})}$.
En effet, on a alors
$$
\mu_{\mathcal C} \bigl( S(a , b) \bigr) \leq \sum_{k = 1}^{n} \mu_{\mathcal C}\bigl( S_k \bigr)
$$
avec
$$
\mu_{\mathcal C}\bigl( S_k \bigr) \leq
\mu_{D_{r}(\omega_{k})}\bigl( S_k \bigr)
$$
pour tout $1 \leq k \leq n$ en vertu de la proposition~\ref{pptinclusion}~(iv).

Or, chaque rectangle $S_k$ étant la réunion des
adhérences de deux triangles, on a $\mu_{D_{r}(\omega_{k})}(S_k) \leq 2 \pi$
puisque tout triangle est contenu dans un triangle
idéal et que $\pi$ est l'aire d'un triangle idéal dans le disque
hyperbolique. Ceci achève la preuve du lemme~\ref{Recouvrement}.
\end{itemize}
\end{proof}

\begin{address}
  Bruno Colbois \\
  Institut de mathématique 
  Université de Neuchâtel 
  Rue Émile Argand 11
  CH--2007 Neuchâtel 
  Switzerland 
  \email{Bruno.Colbois@unine.ch}
\end{address}

\begin{address}
  Constantin Vernicos \\
  Institut de mathématique 
  Université de Neuchâtel 
  Rue Émile Argand 11
  CH--2007 Neuchâtel 
  Switzerland 
  \email{Constantin.Vernicos@unine.ch}
\end{address}

\begin{address}
  Patrick Verovic \\
  Université de Savoie 
  Campus scientifique 
  Laboratoire de mathématique
  F--73376 Le Bourget-du-Lac cedex 
  France
  \email{Patrick.Verovic@univ-savoie.fr}
\end{address}

\end{document}